\documentclass[11pt]{amsart}
\usepackage{geometry}                
\usepackage{graphicx}
\usepackage{amssymb}
\usepackage{epstopdf}
\usepackage{units}
\usepackage{wasysym}
\newcommand{\beq}{
\begin{eqnarray}
}
\newcommand{\eeq}{
\end{eqnarray}
}

\theoremstyle{plain}

\newtheorem{definition}{Definition}

\newtheorem{exercise}{Exercise}

\newtheorem{proposition}{Proposition}

\title[Toroidal Mosaic Numbers]{On Upper Bounds for Toroidal Mosaic Numbers}
\author{Michael Carlisle}
\address{Baruch College, City University of New York (CUNY)\\
New York, NY \ 10010 \ \ USA}
\email{michael.carlisle@baruch.cuny.edu}
\author{Michael S. Laufer}
\address{No affiliation.}
\email{michaelswanlaufer@gmail.com}
\begin{document}
\maketitle

\begin{abstract}
In this paper, we construct mosaic representations of knots on the torus, rather than in the plane. This consists of a particular choice of the ambient group $\mathbb{A}$, as well as different definitions of \emph{contiguous} and \emph{suitably connected}. We present conditions under which mosaic numbers might decrease by this projection, and present a tool (called \emph{waste}) to measure this reduction. 
We show that the order of edge identification in construction of the torus sometimes yields different resultant knots from a given mosaic when reversed. Additionally, in the Appendix we give the catalog of all torus $2$-mosaics. 
\end{abstract}

\section{Introduction}

This paper was inspired from open question (8) of \cite{LK}; unless otherwise noted, definitions come from \cite{LK}. For simplicity of exposition, we will frequently use the term ``knot" to mean either a knot or a link, and we adopt the conventions put forth in \cite{LK}, making the following adjustments to definitions, and using hats $\,\widehat{}\,$ to distinguish toroidal objects from planar objects:

\begin{definition}
A \textbf{toroidal $n$-mosaic} is an $n$-mosaic projected onto the torus in $\mathbb{R}^3$ by identifying opposite edges of the $n$-mosaic. (We examine these mosaics via their matrix representation as given in \cite{LK}.)
\end{definition}

\begin{definition}
Two tiles in a toroidal $n$-mosaic are said to be \textbf{contiguous} if they lie immediately next to each other in either the same row or the same column, or are on opposite ends of a row or column (i.e. the tiles in rows 0 and $n-1$ in column $j$ are contiguous, as are the tiles in columns 0 and $n-1$ of row $j$). An unoriented tile within a mosaic is said to be \textbf{toriodally suitably connected} if each of its connection points touches a connection point of a contiguous tile.
\end{definition}

\begin{definition}
A \textbf{toroidal knot $n$-mosaic} is a toroidal $n$-mosaic in which every tile is toroidally suitably connected. The set of toroidal knot $n$-mosaics is denoted $\widehat{\mathbb{K}^{(n)}}$.
\end{definition}
\begin{definition}
The \textbf{toroidal knot mosaic ambient group} $\widehat{\mathbb{A}(n)}$ is the group of all permutations of $\widehat{\mathbb{K}^{(n)}}$ generated by mosaic planar isotopy moves, mosaic Reidemeister moves, and cyclic permutation of rows and columns.
\end{definition}

\begin{figure}[h]
   \centering
\[
\begin{tabular}
[c]{c}%
$%
\begin{array}
[c]{cc}%
{\includegraphics[
height=0.3269in,
width=0.3269in
]%
{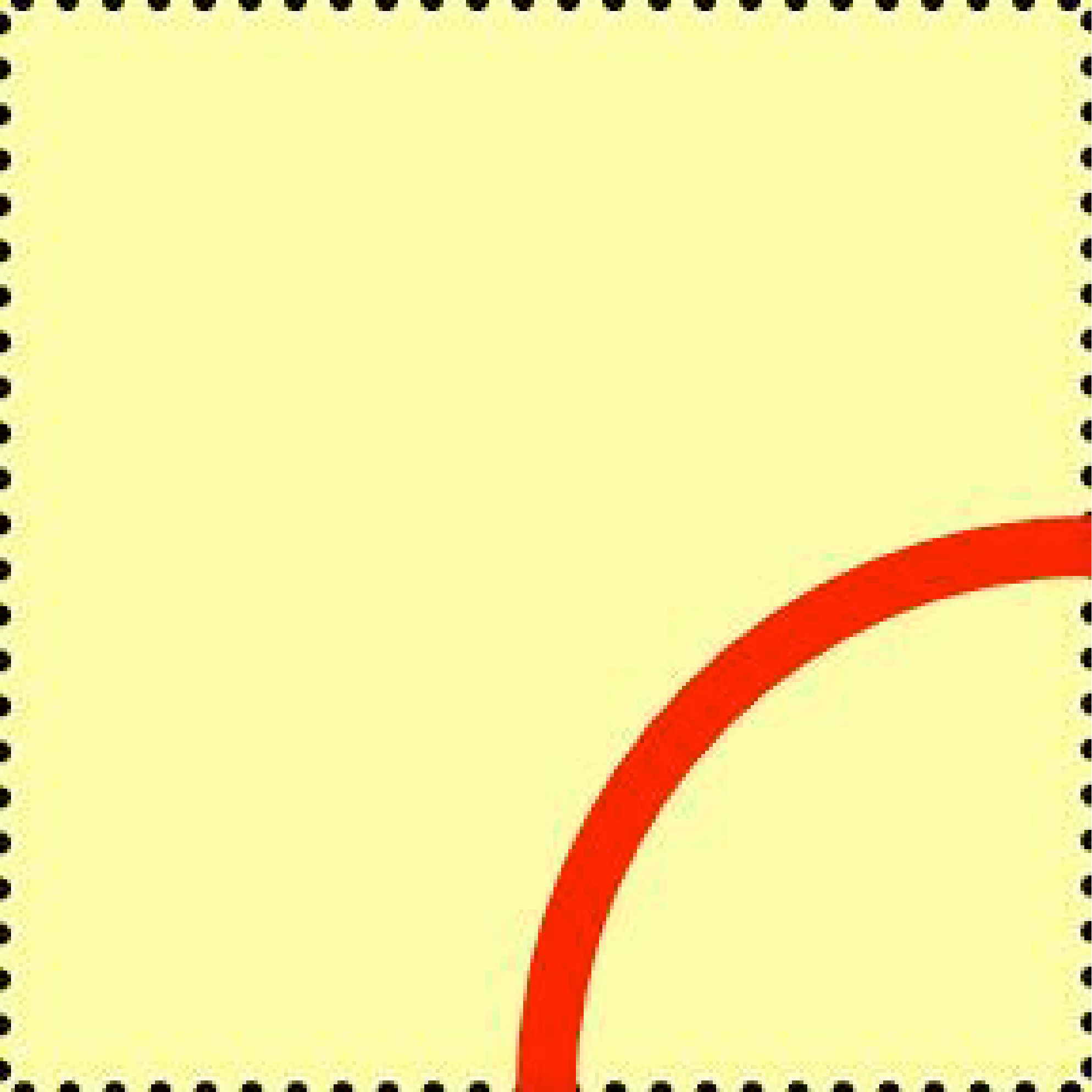}%
}%
&
{\includegraphics[
height=0.3269in,
width=0.3269in
]%
{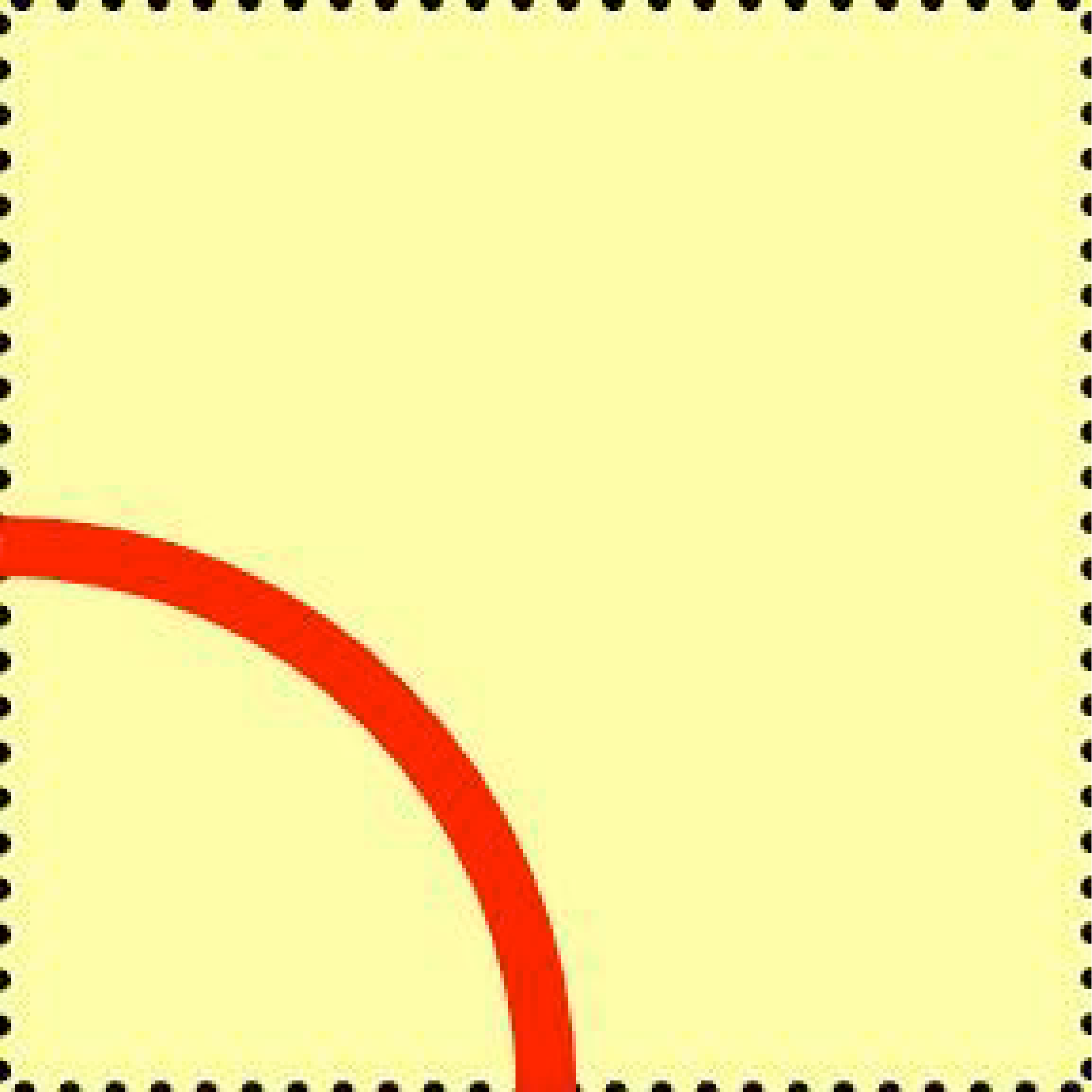}%
}%
\\%
{\includegraphics[
height=0.3269in,
width=0.3269in
]%
{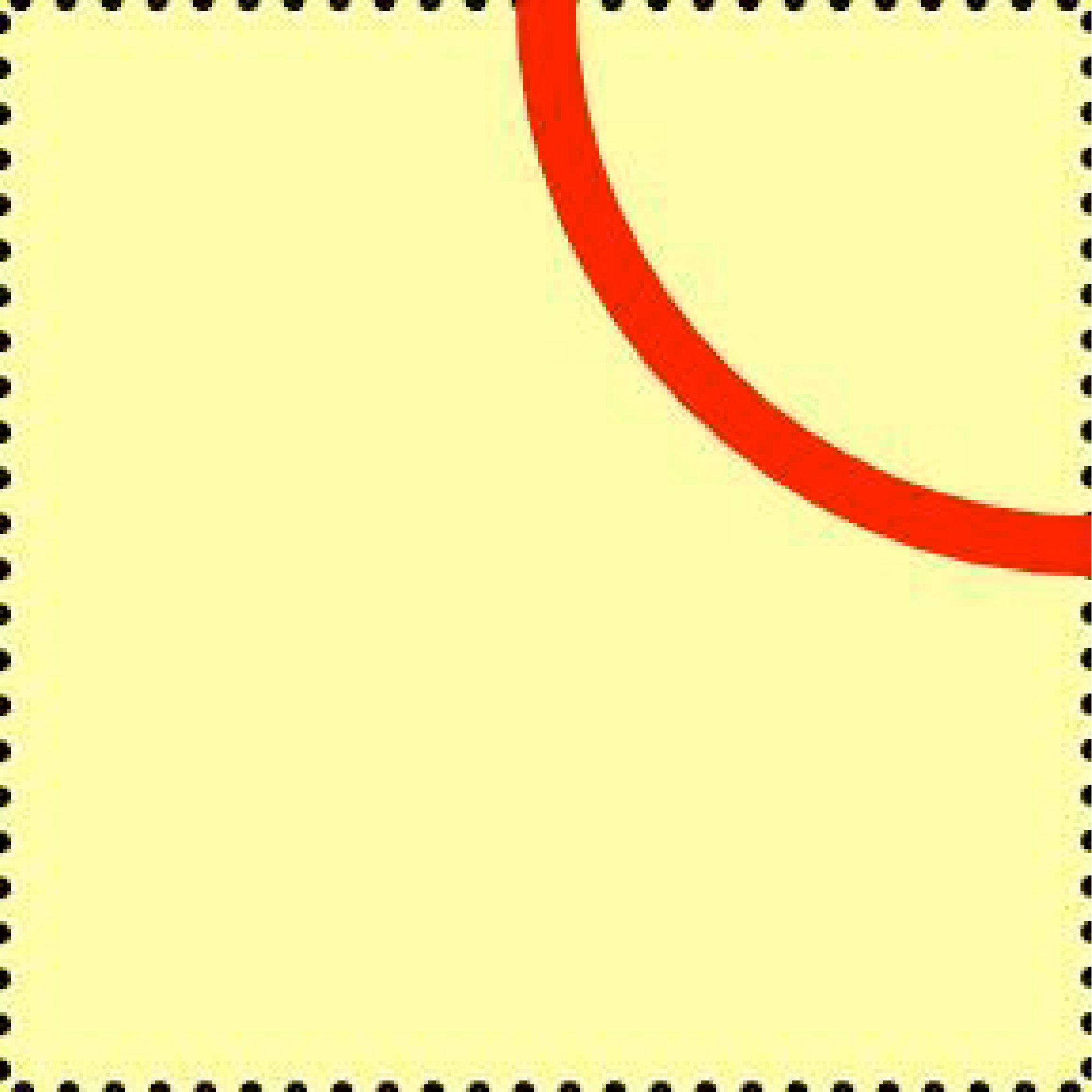}%
}%
&
{\includegraphics[
height=0.3269in,
width=0.3269in
]%
{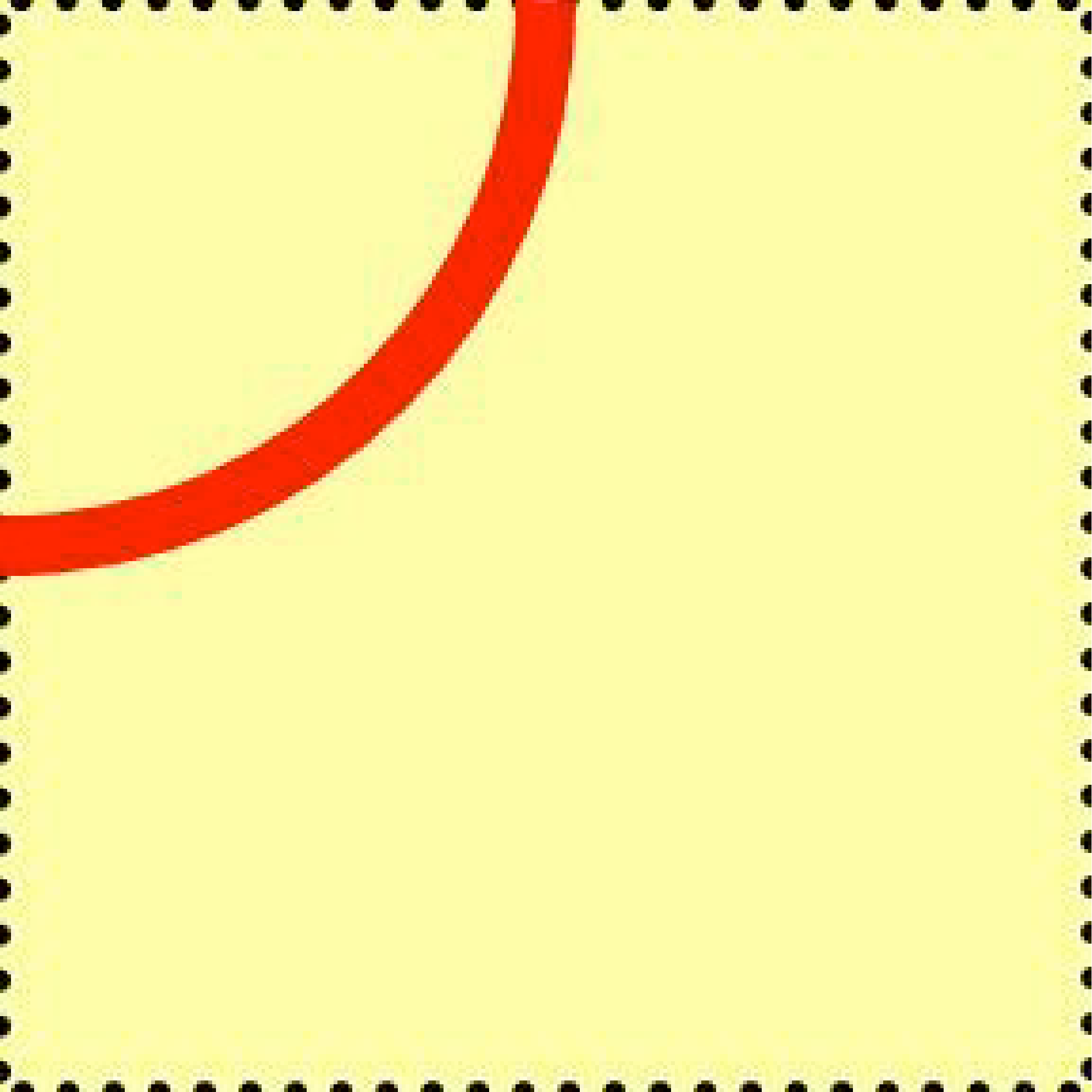}%
}%
\end{array}
\longleftrightarrow%
\begin{array}
[c]{cc}%
{\includegraphics[
height=0.3269in,
width=0.3269in
]%
{ut01.pdf}%
}%
&
{\includegraphics[
height=0.3269in,
width=0.3269in
]%
{ut02.pdf}%
}%
\\%
{\includegraphics[
height=0.3269in,
width=0.3269in
]%
{ut04.pdf}%
}%
&
{\includegraphics[
height=0.3269in,
width=0.3269in
]%
{ut03.pdf}%
}%
\end{array}
$%
\end{tabular}
\]
   \caption{$K_{24} \, \longleftrightarrow \, K_{24} B_2$}
 \label{fig:planartoral}
\end{figure}

``Cyclic permutations of rows and columns'' refers to matrix multiplication of the $n$-mosaic $M \in \mathbb{M}^{(n)}$ by the $n \times n$ (unitary) permutation matrix $B_n$ (on the left for row shifts, on the right for column shifts): 
\[ B_n = \begin{bmatrix}
0 & 1 & 0 & \cdots & 0 \\
0 & 0 & 1 & \cdots & 0 \\
\vdots & \vdots & \vdots & \ddots & \vdots \\
0 & 0 & 0 & \cdots & 1 \\
1 & 0 & 0 & \cdots & 0 \\
\end{bmatrix}. \]
The addition of cyclic permutations of planar knot $n$-mosaics implies that toroidal knot $n$-mosaic set $\widehat{\mathbb{K}^{(n)}}$ contains elements of $\mathbb{M}^{(n)}$ that have planarly non-contiguous tiles; i.e. $\mathbb{K}^{(n)} \subset \widehat{\mathbb{K}^{(n)}} \subset \mathbb{M}^{(n)}$. See Figure \ref{fig:planartoral} for an example: $K_{24} \in \mathbb{K}^{(2)}$, but $K_{24} B_2 \in \widehat{\mathbb{K}^{(2)}} \setminus \mathbb{K}^{(2)}$.

We now state a variation of the notion of mosaic number given in \cite{LK}, open question (8), under this toroidal paradigm: 

\begin{definition}
The \textbf{toroidal mosaic number} of a knot $k$ is the smallest integer $n$ such that $k$ is representable as a toroidal knot $n$-mosaic. 
\end{definition}

Note that toroidal knot $n$-mosaics are two-dimensional projections of three-dimensional knots. If we instead make this projection onto mosaic tiles on a torus (which itself is a three-dimensional object representable in two dimensions), we can lower this mosaic number.

\section{Waste, Density, Embedding}

The difference in planar and toroidal mosaic number can be captured succinctly in the concept of \emph{waste}. 
\begin{definition}
The \textbf{(normalized) waste} of a mosaic tile is $\nicefrac{1}{4}$ the number of tile edges without connection points. The \textbf{total waste} of a knot $n$-mosaic is the sum of its tiles' waste.
\end{definition}
The blank tile $T_0$ has waste 1, the line and 1-corner tiles $T_1$ through $T_6$ have waste $\nicefrac{1}{2}$; crossings and 2-corner tiles ($T_7$ through $T_{10}$) have waste 0. 

\begin{definition}
A knot $n$-mosaic is called \textbf{dense} if it has total waste 0. 
\end{definition}

\begin{figure}[h]
   \centering
   \begin{eqnarray}
\begin{array}[c]{cccc}
{\includegraphics[
height=0.3269in,
width=0.3269in
]%
{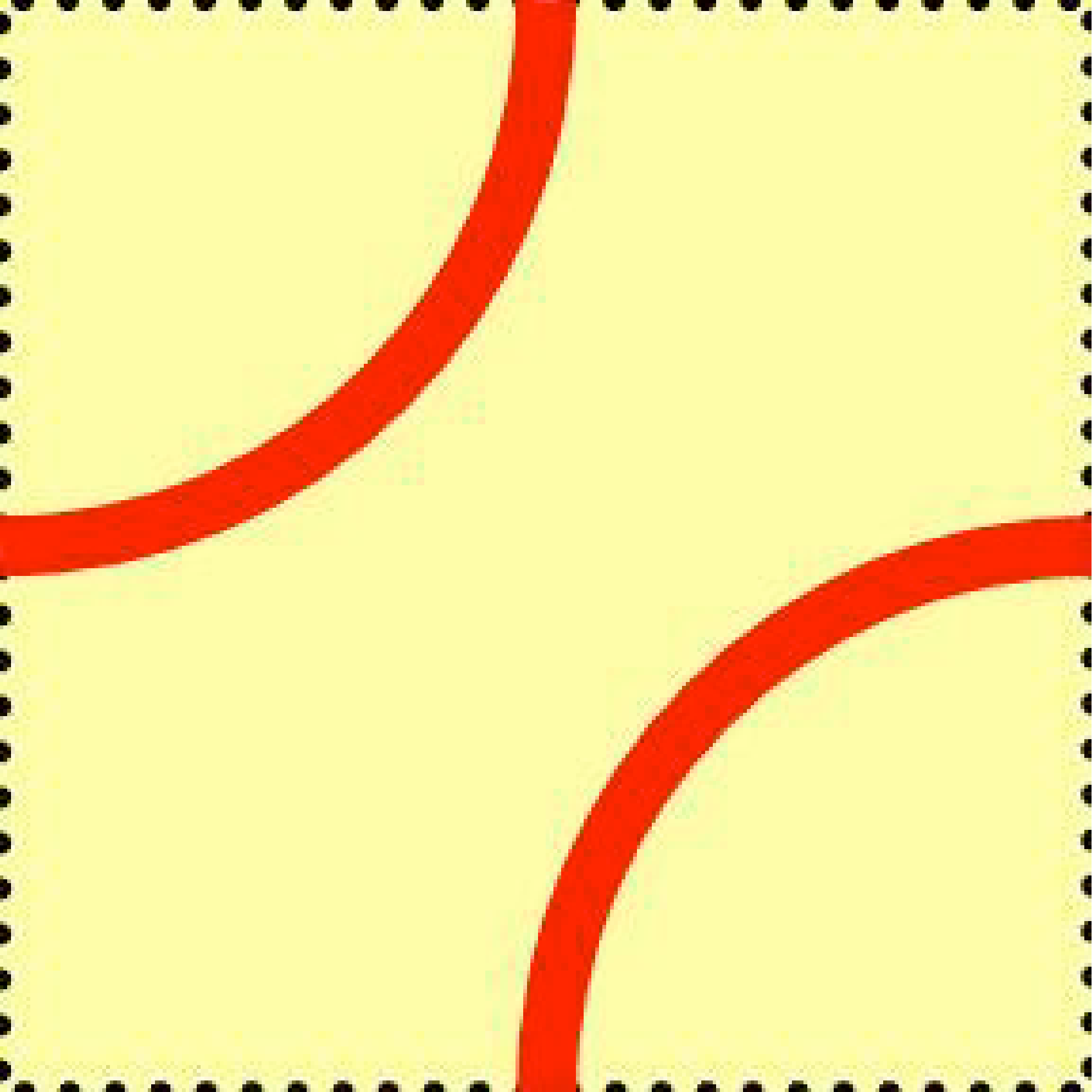}%
}%
& {\includegraphics[
height=0.3269in,
width=0.3269in
]%
{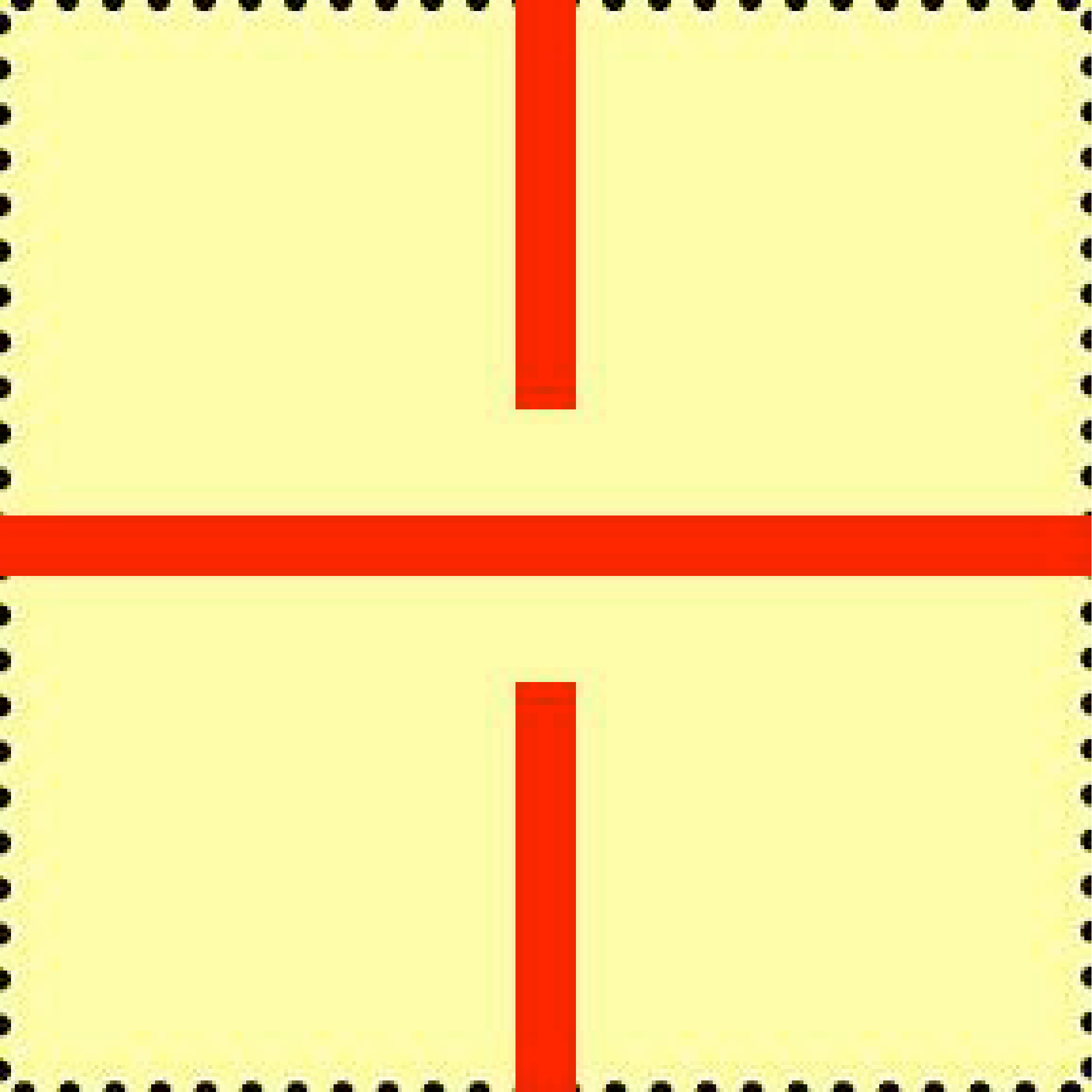}%
}%
& {\includegraphics[
height=0.3269in,
width=0.3269in
]%
{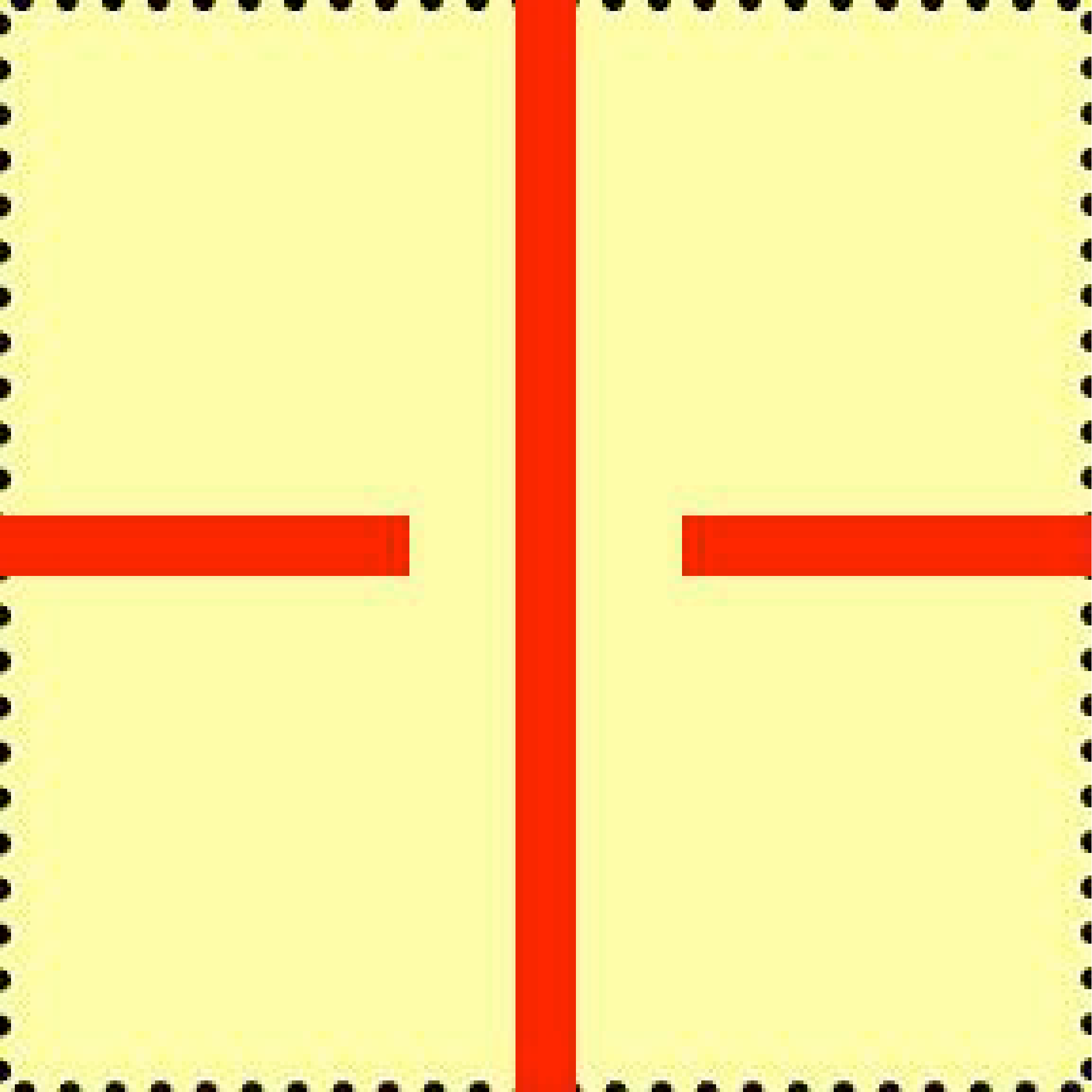}%
}%
 \\
{\includegraphics[
height=0.3269in,
width=0.3269in
]%
{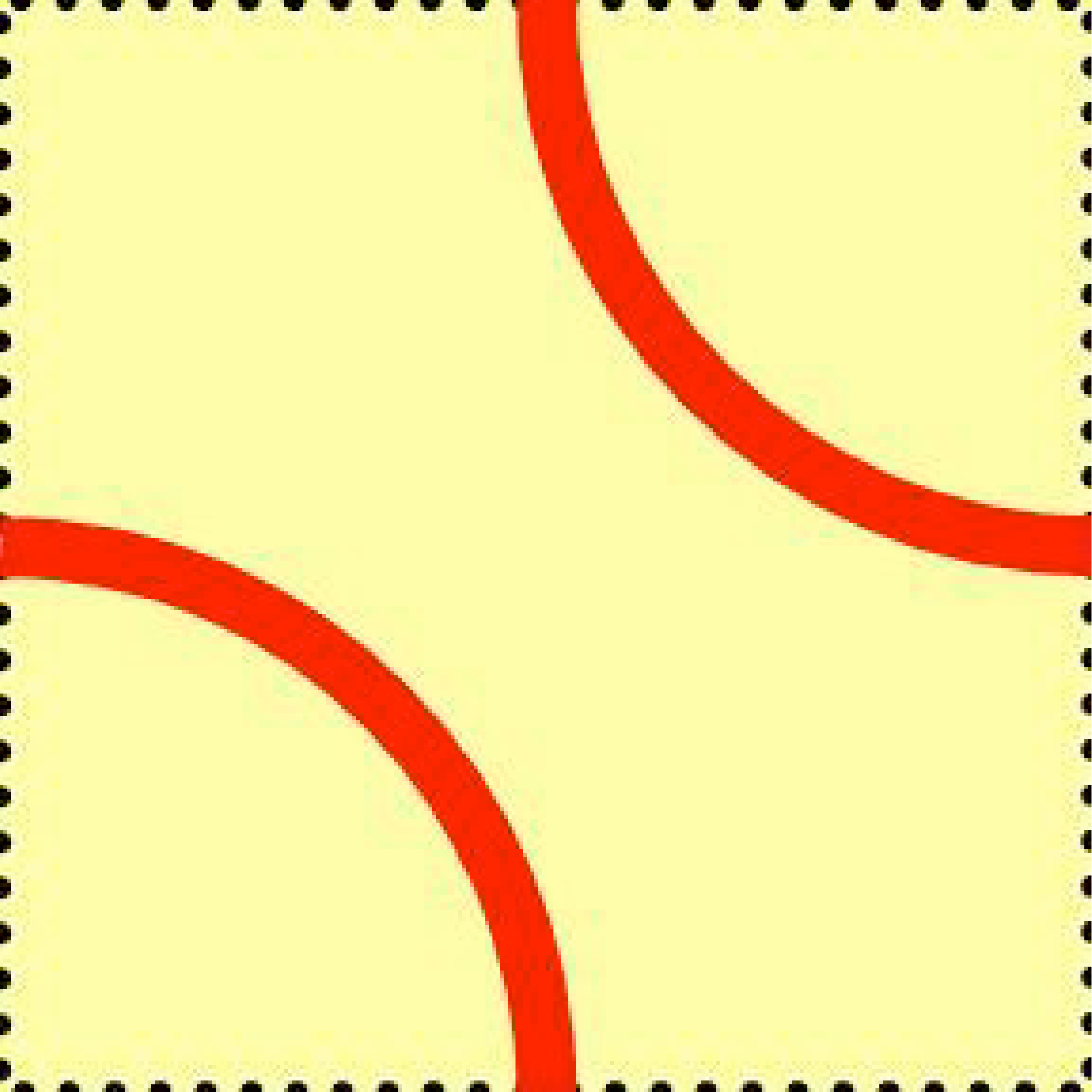}%
}%
& {\includegraphics[
height=0.3269in,
width=0.3269in
]%
{ut07.pdf}%
}%
& {\includegraphics[
height=0.3269in,
width=0.3269in
]%
{ut09.pdf}%
}%
 \\{\includegraphics[
height=0.3269in,
width=0.3269in
]%
{ut09.pdf}%
}%
& {\includegraphics[
height=0.3269in,
width=0.3269in
]%
{ut07.pdf}%
}%
& {\includegraphics[
height=0.3269in,
width=0.3269in
]%
{ut08.pdf}%
}%
\end{array}\nonumber
\end{eqnarray}
   \caption{A dense presentation of $3_1$ on the torus with 4 crossings.}
   \label{fig:dense1}
\end{figure}

Having waste 0 is not sufficient for a tiling to be a representative example demonstrating the mosaic number of a given link; see Figure \ref{fig:dense1} for an example. Note that planar mosaic knots will never be dense, as each of the tiles on the boundary will have at least waste $\nicefrac{1}{2}$. This implies that in larger presentations of planar knots waste will increase at least linearly. It is suggestive that more waste in the plane means more room for reducing waste by moving to the torus.

Figure \ref{fig:dense1} is an example of the $3_1$ knot on the torus with 4 crossings in the diagram. The lower left crossing is superfluous.

\begin{figure}[h]
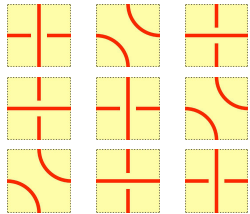

   \centering
   \begin{eqnarray}
\begin{array}[c]{cccc}
{\includegraphics[
height=0.3269in,
width=0.3269in
]%
{}%
}%
& {\includegraphics[
height=0.3269in,
width=0.3269in
]%
{}%
}%
& {\includegraphics[
height=0.3269in,
width=0.3269in
]%
{}%
}%
 \\
{\includegraphics[
height=0.3269in,
width=0.3269in
]%
{}%
}%
& {\includegraphics[
height=0.3269in,
width=0.3269in
]%
{}%
}%
& {\includegraphics[
height=0.3269in,
width=0.3269in
]%
{}%
}%
 \\{\includegraphics[
height=0.3269in,
width=0.3269in
]%
{}%
}%
& {\includegraphics[
height=0.3269in,
width=0.3269in
]%
{}%
}%
& {\includegraphics[
height=0.3269in,
width=0.3269in
]%
{}%
}%
\end{array}\nonumber
\end{eqnarray}
   \caption{A Dense Presentation of the Borromean Rings on the Torus}
   \label{fig:dense2}
\end{figure}

\begin{figure}
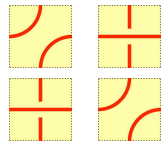

\begin{eqnarray}
\begin{array}[c]{cccc}
{\includegraphics[
height=0.3269in,
width=0.3269in
]%
{}%
}%
& {\includegraphics[
height=0.3269in,
width=0.3269in
]%
{}%
}%
 \\
{\includegraphics[
height=0.3269in,
width=0.3269in
]%
{}%
}%
& 
{\includegraphics[
height=0.3269in,
width=0.3269in
]%
{}%
}%
\end{array}\nonumber
\end{eqnarray}
 \caption{Possible 2-crossing Hopf link on the torus.}
   \label{fig:hopf1}
\end{figure}

\pagebreak

\begin{definition} \label{def:inj}
The \textbf{toroidal mosaic injection} for a toroidal knot $n$-mosaic is defined by  
\begin{align*}
\widehat{\iota}:\mathbb{M}^{(n)} & \longrightarrow \mathbb{M}^{(n+1)} \\
M^{(n)} & \longmapsto M^{(n+1)}  
\end{align*}
as
\end{definition}
\[
M_{ij}^{(n+1)}=\left\{
\begin{array}
[c]{cl}%
M_{ij}^{(n)} & \text{if }0\leq i,j<n\\
& \\%
\raisebox{-0.1003in}{\includegraphics[
height=0.3269in,
width=0.3269in
]%
{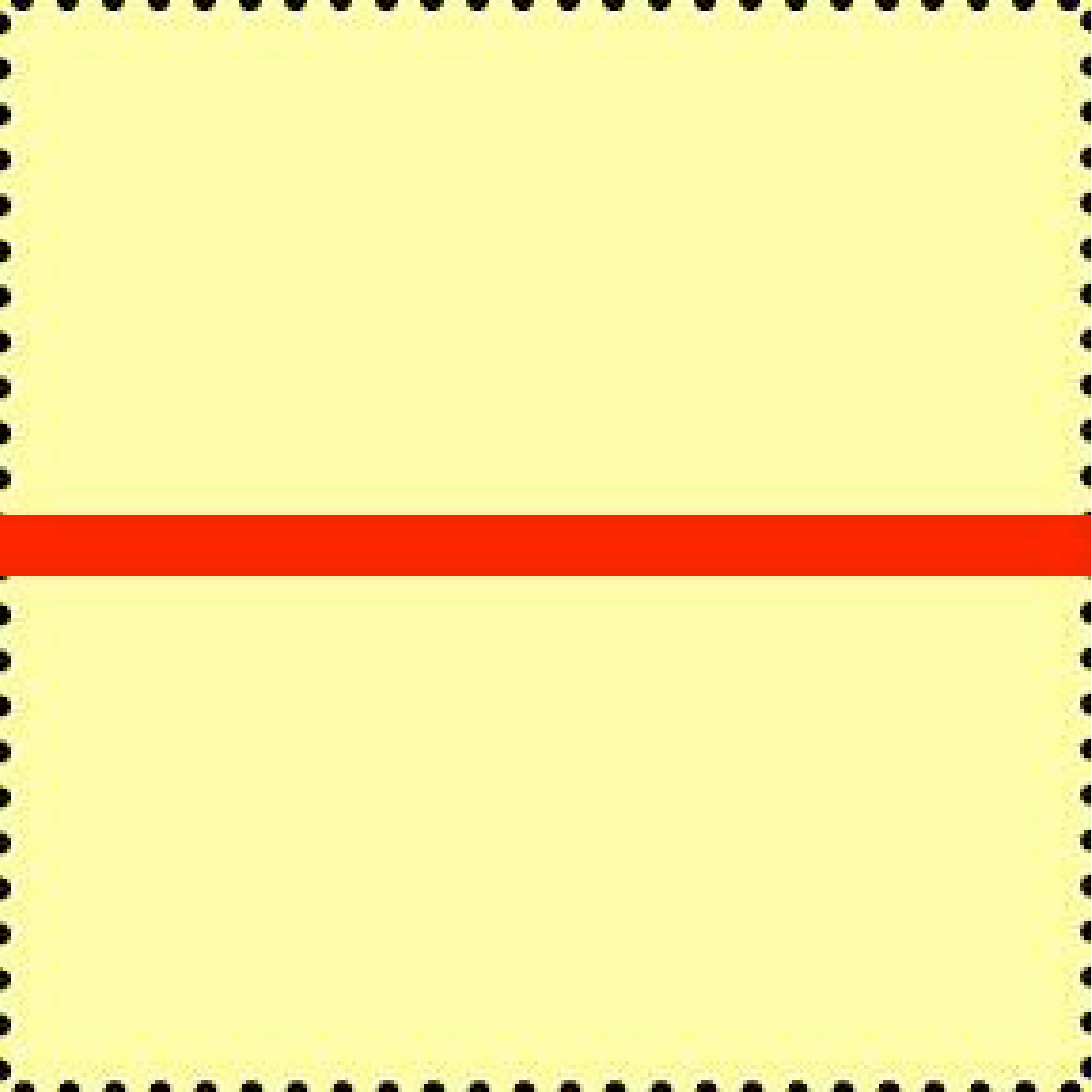}%
}%
& = T_5 \text{ if } j=n, M_{i(n-1)}^{(n)} \text{ has a connection to the right } \\
& \\
\raisebox{-0.1003in}{\includegraphics[
height=0.3269in,
width=0.3269in
]%
{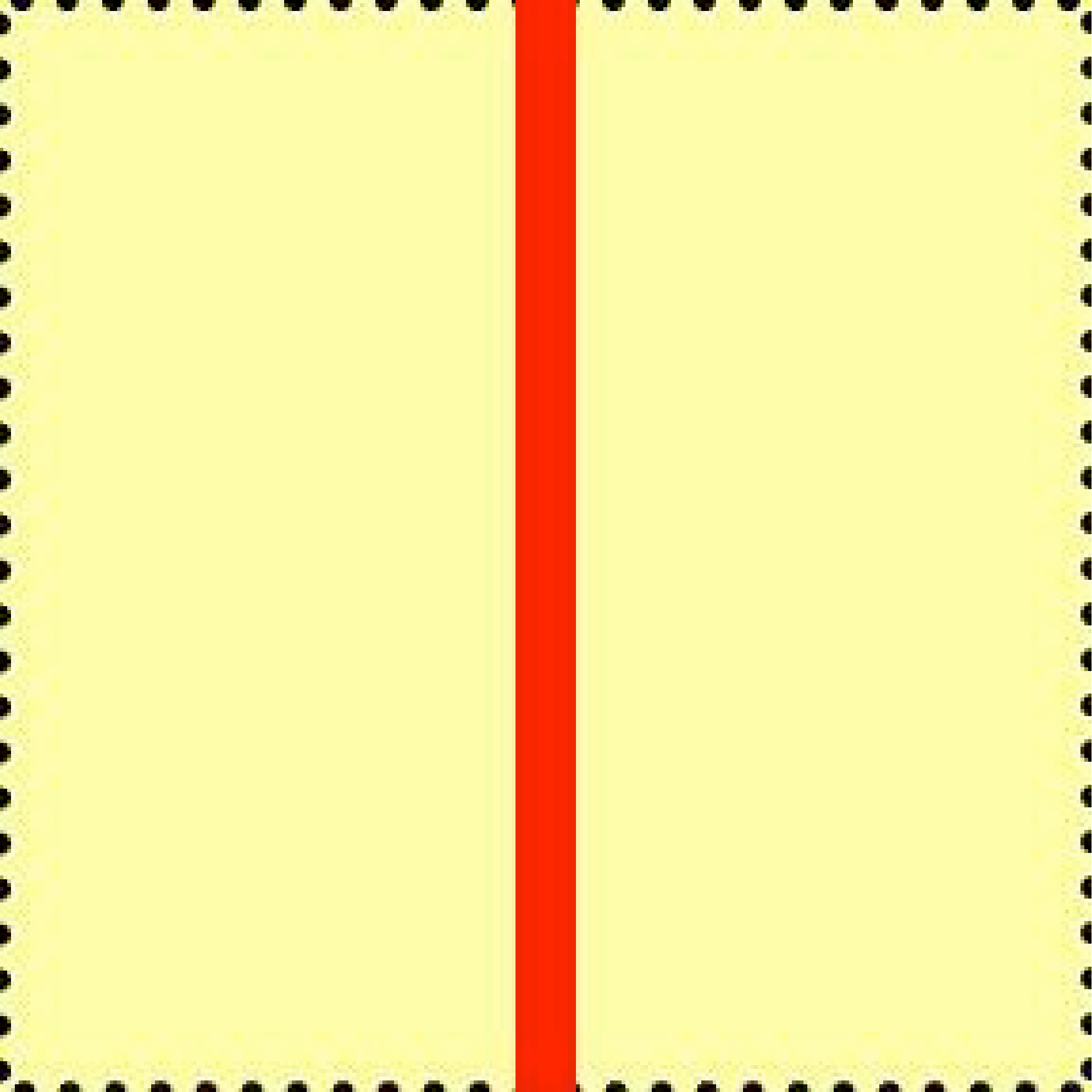}%
}%
& = T_6 \text{ if } i=n, M_{(n-1)j}^{(n)} \text{ has a connection on the bottom } \\%
& \\
\raisebox{-0.1003in}{\includegraphics[
height=0.3269in,
width=0.3269in
]%
{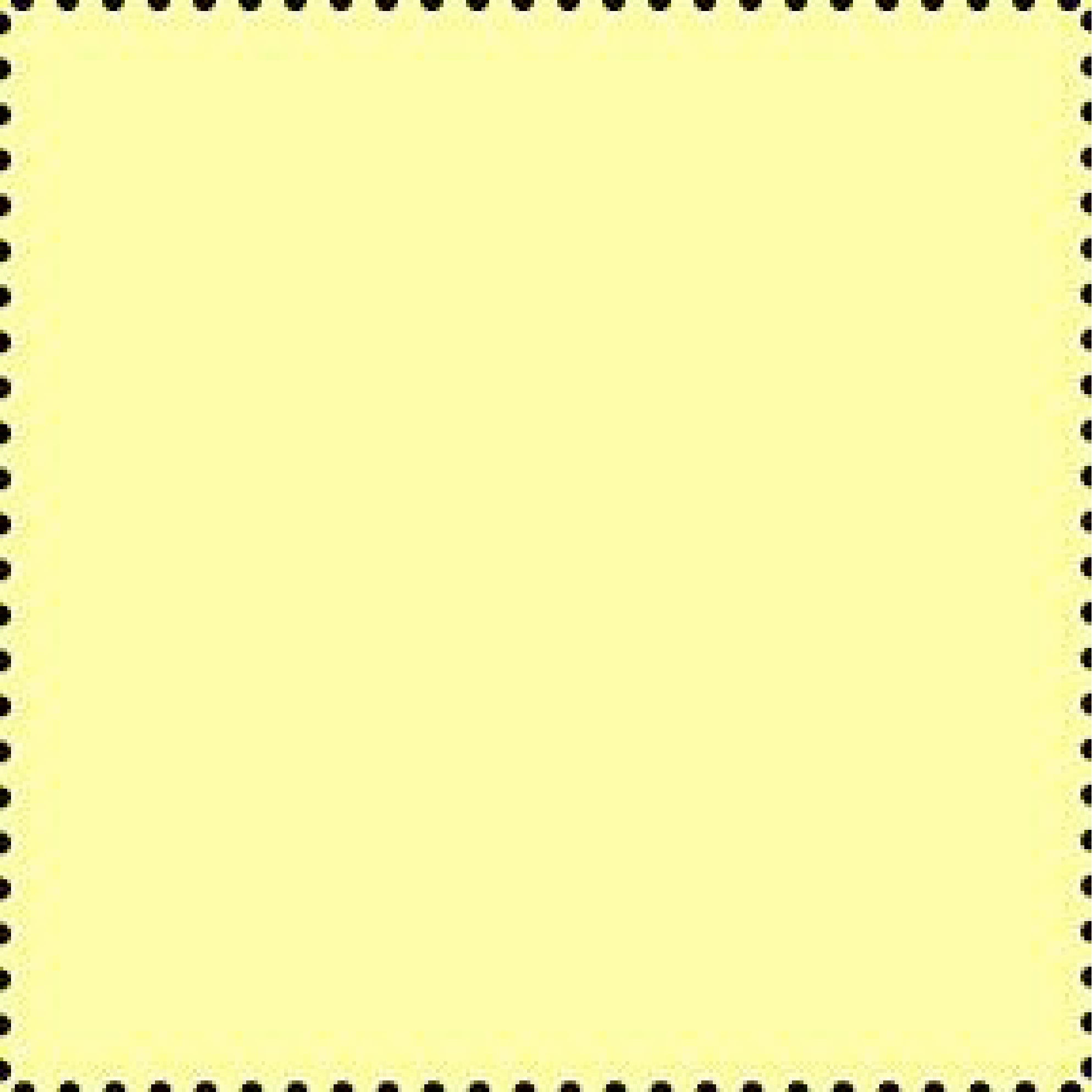}%
}%
& = T_0 \text{ otherwise.}%
\end{array}
\right.
\]

\bigskip

We now give bounds on the amount of waste induced by embedding a toroidal knot $n$-mosaic. 

\begin{proposition}
Let $M \in \widehat{\mathbb{K}^{(n)}}$. Then, if $w(M)$ is the amount of waste of $M$, 
\[ w(M) + n + 1 \leq w(\widehat{\iota}(M)) \leq w(M) + 2n+1, \]
where the upper bound is achieved if $M \in \mathbb{K}^{(n)}$.
\end{proposition}

\noindent \textbf{Proof} A toroidal embedding always uses the blank tile $T_0$ (with waste 1) in cell $(n,n)$. Considering $M$ as an unoriented planar mosaic, the waste upper bound is achieved when every edge on $M$'s boundary is wasted. Since $M \in \widehat{\mathbb{K}^{(n)}}$, this means $M \in \mathbb{K}^{(n)}$, and so the toroidal embedding matches the planar embedding, using a blank tile for every new position (a total of $2n+1$ blank tiles).

Again considering $M$ as an unoriented planar mosaic, if $M$ has no waste on its boundary (which does not necessarily mean it is dense; it may have waste in its interior), then its embedding uses $T_5$ (with waste $\frac{1}{2}$ per tile) for the first $n$ entries in column $n$ and $T_6$ (with waste $\frac{1}{2}$ per tile) for the first $n$ entries in row $n$, and  $T_0$ in cell $(n,n)$. $\qed$

This injection extends graded system $(\mathbb{K},\mathbb{A})$ given in \cite{LK} to the torus. The symbol $\widehat{\mathbb{K}}$ denotes the directed system of sets $\{\widehat{\mathbb{K}^{(n)}}\longrightarrow\widehat{\mathbb{K}^{(n+1)}}:n=1,2,3\cdots\}$ and $\widehat{\mathbb{A}}$ denotes the directed system of groups $\{\widehat{\mathbb{A}(n)}\longrightarrow\widehat{\mathbb{A}(n+1)}:n=1,2,3\cdots\}$, thus
\[
\left(  \widehat{\mathbb{K}},\widehat{\mathbb{A}}\right)
  = \left(\widehat{\mathbb{K}^{(1)}},\widehat{\mathbb{A}(1)}\right) 
\longrightarrow\left(
\widehat{\mathbb{K}^{(2)}},\widehat{\mathbb{A}(2)}\right)
\longrightarrow \cdots
\longrightarrow \left(
\widehat{\mathbb{K}^{(n)}},\widehat{\mathbb{A}(n)}\right)
  \longrightarrow\cdots
\]

\section{A Strange Paradox}

Projecting onto the torus introduces ``hidden'' crossings that may destroy the well-definedness of a knot mosaic. Consider the $1 \times 1$ mosaic of one crossing. On the torus, is it the Hopf link?

\begin{figure}[h]
   \centering
{\includegraphics[
height=0.3269in,
width=0.3269in
]%
{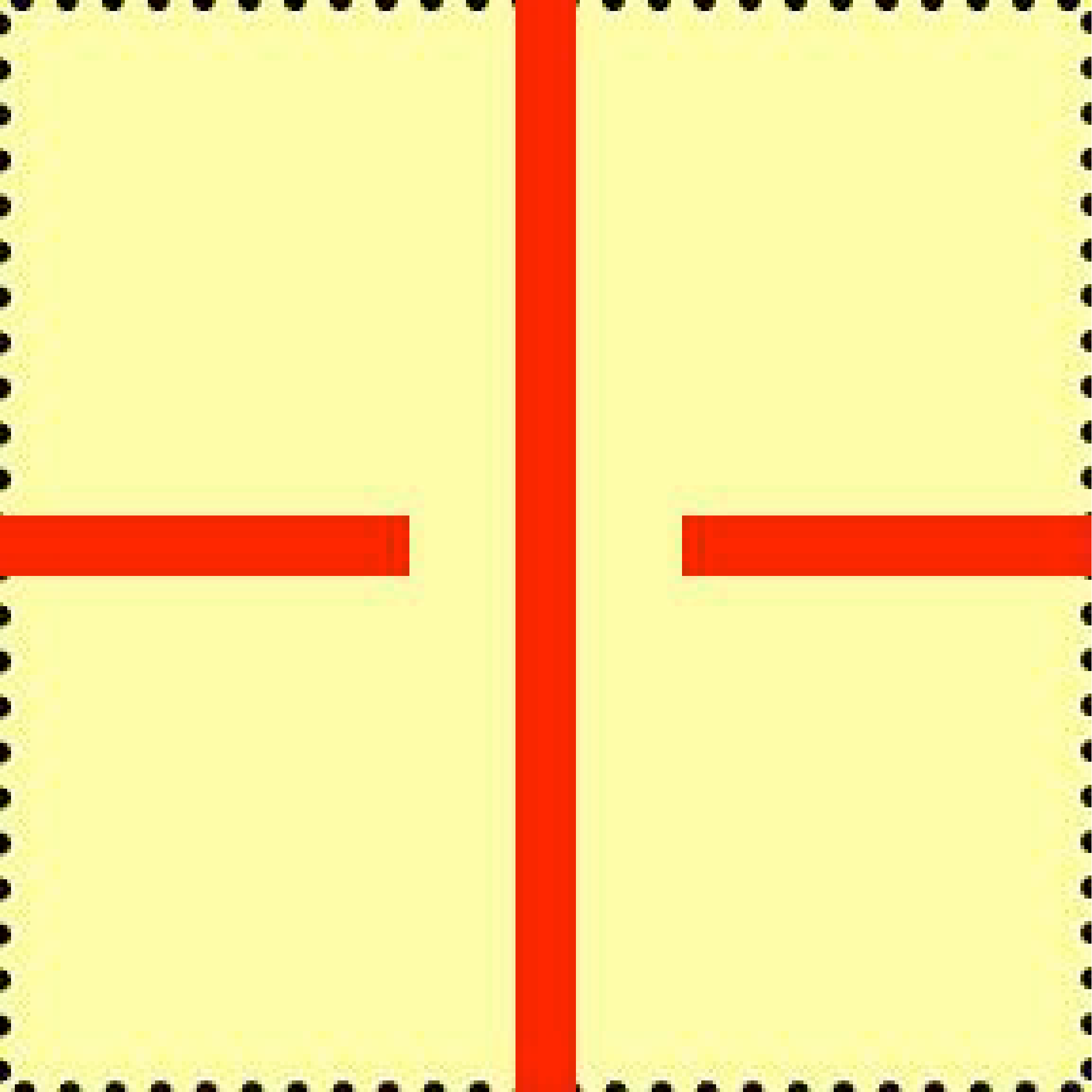}%
}
   \caption{Is this also the Hopf Link on the torus?}
   \label{fig:hopf2}
\end{figure}

The answer is that it depends on which pair of edges you connect first when constructing the torus. Figure \ref{fig:hopf2} will either be the Hopf link, or two disjoint unknots, depending on the choice. Connecting top and bottom first will yield the Hopf link; connecting left and right first will yield two disjoint unknots. This forces us to reconsider our previous constructions, and see if we haven't unknowingly introduced crossings ``off the mosaic". 

It becomes more complicated to see the implications of embedding the knot on the torus in $\mathbb{R}^3$ than to merely identify opposite edges. 

\begin{figure}
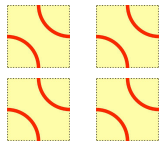

\begin{eqnarray}
\begin{array}[c]{cccc}
{\includegraphics[
height=0.3269in,
width=0.3269in
]%
{}%
}%
& {\includegraphics[
height=0.3269in,
width=0.3269in
]%
{}%
}%
 \\
{\includegraphics[
height=0.3269in,
width=0.3269in
]%
{}%
}%
& 
{\includegraphics[
height=0.3269in,
width=0.3269in
]%
{}%
}%
\end{array}\nonumber
\end{eqnarray}
 \caption{Hopf Link on the 3D torus with no explicit crossings.}
   \label{fig:hopf3}
\end{figure}

Because we introduce crossings by the topology of the torus, it is possible to construct links which do not show any crossings in the diagram. The Hopf link is diagrammed on Figure \ref{fig:hopf3} with no explicit crossings.

\begin{figure}[h]
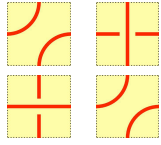

   \centering
\begin{eqnarray}
\begin{array}[c]{cccc}
{\includegraphics[
height=0.3269in,
width=0.3269in
]%
{}%
}%
& {\includegraphics[
height=0.3269in,
width=0.3269in
]%
{}%
}%
 \\
{\includegraphics[
height=0.3269in,
width=0.3269in
]%
{}%
}%
& 
{\includegraphics[
height=0.3269in,
width=0.3269in
]%
{}%
}%
\end{array}\nonumber 
\end{eqnarray} 
  \caption{The actual 2-crossing Hopf link on the torus (differs from Figure \ref{fig:hopf1})}
  \label{fig:hopf4}
\end{figure}

We note that in this case order of edges does not matter, because of the symmetry of this presentation. Figure \ref{fig:hopf4} is the Hopf link when top and bottom edges are identified first, and disjoint unknots when left and right are identified first.

In order to remove this ambiguity one of two possible conventions must be adopted: Convention 1.  Left--to--right is designated as the meridianal direction. (Call these \textbf{meridianal toroidal mosaic knots}.)  
 Convention 2.  Top--to--bottom is designated as the longitudinal direction. (Call these \textbf{longitudinal toroidal mosaic knots}.)  

Henceforth, we shall implicitly adopt the latter convention.

Note that the catalog of $n$-mosaics changes markedly on the torus because cyclic translation is an equivalence relation. For instance, in the catalog put forth in \cite{LK}, the diagrams labeled as $K_1, K_2, K_4$ and $K_{11}$ would all be toroidally equivalent. However, the toroidal representation of the Borromean rings in Figure \ref{fig:dense2} would not be an allowable figure in the planar 3-mosaic catalog. Hence, the toroidal $n$-mosaic catalog, surprisingly, is much larger than the planar $n$-mosaic catalog.

Rotations of $n$-mosaics, however, as we have seen, yield very different results, so each mosaic must be rotated by $\frac{\pi}{2}$ and reexamined to see if it has differing structure.

The mosaics which are the smallest representations of knots and links on the torus thus far all happen to be dense mosaics. However, density as a mosaic property does not seem to imply any other properties at this point.

\begin{proposition}
The dense toroidal 1-mosaic $K_5$ is the smallest possible toroidal mosaic presentation of the Hopf link.
\end{proposition}

\begin{proposition}
The dense toroidal 2-mosaics $K_{49}$ and $K_{53}$ are the smallest possible toroidal mosaic presentations of $3_1$.
\end{proposition}

\begin{proposition}
The dense toroidal 2-mosaics $K_{73}, K_{75}, K_{83}, K_{85}$ and $K_{88}$ are the smallest possible toroidal mosaic presentations of the $4_1^2$ link.
\end{proposition}

\begin{proposition}
The dense toroidal 2-mosaics $K_{93}$ and $K_{94}$ are the smallest possible toroidal mosaic presentations of the $6_1^3$ link.
\end{proposition}

\begin{proposition}
The dense toroidal 2-mosaic $K_{98}$ is the smallest possible toroidal mosaic presentation of the $8_4^3$ link.
\end{proposition}

To see that these propositions are true, one can examine the catalog. \\ \\
Each possible toroidal $2$-mosaic was constructed, examined to determine which knot it represented, and equivalent diagrams were removed. Thus far we have not uncovered any simple methodology to determine what link a diagram represents, or if two diagrams are equivalent.
\\ \\
Ironically, the so-called ``Torus Knots" which come from wrapping a closed path around the torus moving at particular rates meridianally and longitudinally, all require very large mosaics to construct as toroidal mosaics.

\section{Quantum Toroidal Mosiac Knots}

The addition of projecting knot mosaics onto the torus does not interfere with the construction of a quantum knot system in an analogous way to what was done in section 3.1 of \cite{LK}. Using the same definition of the Hilbert space $\widehat{\mathcal{M}^{(n)}} = \mathcal{M}^{(n)}$ of $n$-mosaics, and the same 11 dimensional Hilbert space, we end up with the same induced basis.

We define the \textbf{Hilbert space of toroidal knot }$\mathbf{n}%
$\textbf{-mosaics }$\widehat{\mathcal{K}^{(n)}}$ as the sub-Hilbert space $\widehat{\mathcal{K}^{(n)}}$ of $\mathcal{M}^{(n)}$ spanned by all orthonormal basis elements labeled by toroidal knot $n$-mosaics.

Toroidal quantum knot systems differ only slightly from their planar counterparts. As noted previously, there are the additional unitary operators which allow for cyclic permutation of rows and columns. Hence, the planar ambient group is entirely contained in the toroidal ambient group, i.e. $\mathbb{A}(n) \subset \widehat{\mathbb{A}(n)}$. 
\\ \\
We define a quantum toroidal knot system, much as in the planar case: 

\begin{definition}
Let $n$ be a positive integer. A \textbf{quantum toroidal knot system} $\widehat{Q}\left(
\widehat{\mathcal{K}^{(n)}},\widehat{\mathbb{A}(n)}\right)$ \textbf{of order }$n$ is a quantum
system with the Hilbert space $\widehat{\mathcal{K}^{(n)}}$ of toroidal knot $n$-mosaics as its
state space, and having the ambient group $\widehat{\mathbb{A}(n)}$ as an accessible
unitary control group. \ The states of the quantum system $\widehat{Q}\left(
\widehat{\mathcal{K}^{(n)}},\widehat{\mathbb{A}(n)}\right)$ are called \textbf{quantum toroidal knots of
order }$n$, and the elements of the ambient group $\widehat{\mathbb{A}(n)}$ are called
\textbf{unitary toroidal knot moves}. \ Moreover, the quantum knot system $\widehat{Q}\left(
\widehat{\mathcal{K}^{(n)}},\widehat{\mathbb{A}(n)}\right)$ of\textbf{ }order\textbf{ }$n$ is a
subsystem of the quantum toroidal knot system $\widehat{Q}\left(
\widehat{\mathcal{K}^{(n+1)}},\widehat{\mathbb{A}(n+1)}\right)$ of order\textbf{ }$n+1$. \ Thus, the quantum toroidal knot
systems $\widehat{Q}\left(\widehat{\mathcal{K}^{(n)}},\widehat{\mathbb{A}(n)}\right)$ collectively become a \textbf{nested} \textbf{sequence of quantum toroidal knot systems} which we will denote simply by $\widehat{Q}\left(  \widehat{\mathcal{K}},\widehat{\mathbb{A}}\right)  $. \ In other words,
\[
\widehat{Q}\left(  \widehat{\mathcal{K}},\widehat{\mathbb{A}}\right)
  = \widehat{Q}\left(\widehat{\mathcal{K}^{(1)}},\widehat{\mathbb{A}(1)}\right) 
\longrightarrow \widehat{Q}\left(
\widehat{\mathcal{K}^{(2)}},\widehat{\mathbb{A}(2)}\right)
\longrightarrow \cdots
\longrightarrow \widehat{Q}\left(
\widehat{\mathcal{K}^{(n)}},\widehat{\mathbb{A}(n)}\right)
  \longrightarrow\cdots
\]

\end{definition}

The nested sequence of toroidal knot systems remains consistent, using the convention of toroidal knot injection from Definition \ref{def:inj}. 

Hamiltonians can be constructed identically to the way they are constructed in \cite{LK}, and quantum toroidal knot invariants can be also be found from observables in the same way as was outlined in \cite{LK}.

The set of quantum knot invariants, as defined in \cite{LK} naturally becomes smaller in the toroidal case, because the toroidal ambient group $\widehat{\mathbb{A}(n)}$ is larger. There will be fewer observables which are invariant under conjugation by all the elements of  the ambient group.  Please note that we use the definition for quantum knot invariant used in \cite{LK} rather than the one typically used in quantum topology.

Having added cyclic translations to the ambient knot group does not interfere with being able to construct the quantum knot system. However, this equivalence makes the process of determining if two toroidal mosaics are of the same \emph{knot type} more difficult. 

$\widehat{\mathcal{K}^{(1)}}$ is 7 dimensional, as can be seen in Appendix A, and $\widehat{\mathcal{K}^{(2)}}$ is 97 dimensional. This is quite an increase from the planar cases, which are only 1 and 2 dimensional respectively \cite{LK}. Even $\mathcal{K}^{(3)}$ in the planar case is only 22 dimensional \cite{LK}.

Problem 9 from \cite{LK}, repackaged for torus knots, can be stated:
 
\begin{exercise}
Let $\widehat{D_{n}}$ denote the dimension of the Hilbert space $\widehat{\mathcal{K}^{(n)}}$ of toroidal quantum knot $n$-mosaics. We have shown that $\widehat{D_{1}}=7$ and $\widehat{D_{2}}=97$.  It would be interesting to find $\widehat{D_{n}}$ for other values of $n$.
\ \ A very loose upper bound for $\widehat{D_{n}}$ is obviously $11^{n^{2}}$.\bigskip
\end{exercise}

\section{Conclusions and Further Work}
Moving the model of mosaic knots onto the torus changes the size of the space in which we are working, and allows the theory to remain functional. The fact that edge identification order matters results in having to include rotations in the catalog of all mosaics, even though cyclic translations have been removed.

The introduction of topological artifacts by edge identification begs the question as to what might occur when projecting onto the sphere, the Klein bottle, or the projective plane. These possibilites are under investigation.

It is also intriguing to note that grid diagrams for Floer homology calculations are also discrete grids that produce knot diagrams under an ordered grid connection paradigm (vertical over horizontal), and are invariant under cyclic row and column shifts (i.e. are knot diagrams on the torus). Connections between these two topics are worthy of further research.

\bigskip

\section*{Appendix A: The catalog of all knot 1- and 2-mosaics on the torus.}

Recall that the convention we have adopted identifies top and bottom edges first (i.e. longitudinally). $2_1^2$ or (L2a1) denotes the Hopf link, $+$ denotes disjoint union, and $\#$ denotes knot sum. Only one representation for each mosaic under row or column shift equivalences is displayed. \\ \\

$\underset{K_{0}=\text{$\phi$}}{\fbox{$%
\begin{array}
[c]{ccc}
{\includegraphics[
height=0.3269in,
width=0.3269in
]%
{ut00.pdf}%
}
\end{array}
$}}
\quad\underset{K_{1}=\text{$0_1$}}{\fbox{$%
\begin{array}
[c]{ccc}
{\includegraphics[
height=0.3269in,
width=0.3269in
]%
{ut05.pdf}%
}
\end{array}
$}}\quad\underset{K_{2}=\text{$0_1$}}{\fbox{$%
\begin{array}
[c]{ccc}
{\includegraphics[
height=0.3269in,
width=0.3269in
]%
{ut06.pdf}%
}
\end{array}
$}}\quad\underset{K_{3}=\text{$0_1$}}{\fbox{$%
\begin{array}
[c]{ccc}
{\includegraphics[
height=0.3269in,
width=0.3269in
]%
{ut07.pdf}%
}
\end{array}
$}}\quad\underset{K_{4}=\text{$0_1$}}{\fbox{$%
\begin{array}
[c]{ccc}
{\includegraphics[
height=0.3269in,
width=0.3269in
]%
{ut08.pdf}%
}
\end{array}
$}}\quad\underset{K_{5}=\text{$2_1^2$ (L2a1)}}{\fbox{$%
\begin{array}
[c]{ccc}
{\includegraphics[
height=0.3269in,
width=0.3269in
]%
{ut10.pdf}%
}
\end{array}
$}}\quad\underset{K_{6}=\text{$0_1+0_1$}}{\fbox{$%
\begin{array}
[c]{ccc}
{\includegraphics[
height=0.3269in,
width=0.3269in
]%
{ut09.pdf}%
}
\end{array}
$}}
$
\\ \\ \\
$\underset{K_{7}=\text{$\phi$}}{\fbox{$%
\begin{array}
[c]{ccc}
{\includegraphics[
height=0.3269in,
width=0.3269in
]%
{ut00.pdf}%
}%
& {\includegraphics[
height=0.3269in,
width=0.3269in
]%
{ut00.pdf}%
}%
 \\
{\includegraphics[
height=0.3269in,
width=0.3269in
]%
{ut00.pdf}%
}%
& 
{\includegraphics[
height=0.3269in,
width=0.3269in
]%
{ut00.pdf}%
}%
\end{array}
$}}
\quad\underset{K_{8}=\text{$0_1$}}{\fbox{$%
\begin{array}
[c]{ccc}
{\includegraphics[
height=0.3269in,
width=0.3269in
]%
{ut05.pdf}%
}%
& {\includegraphics[
height=0.3269in,
width=0.3269in
]%
{ut05.pdf}%
}%
 \\
{\includegraphics[
height=0.3269in,
width=0.3269in
]%
{ut00.pdf}%
}%
& 
{\includegraphics[
height=0.3269in,
width=0.3269in
]%
{ut00.pdf}%
}%
\end{array}
$}}
\quad\underset{K_{9}=\text{$0_1$}}{\fbox{$%
\begin{array}
[c]{ccc}
{\includegraphics[
height=0.3269in,
width=0.3269in
]%
{ut00.pdf}%
}%
& {\includegraphics[
height=0.3269in,
width=0.3269in
]%
{ut06.pdf}%
}%
 \\
{\includegraphics[
height=0.3269in,
width=0.3269in
]%
{ut00.pdf}%
}%
& 
{\includegraphics[
height=0.3269in,
width=0.3269in
]%
{ut06.pdf}%
}%
\end{array}
$}}
\quad\underset{K_{10}=\text{$0_1+0_1$}}{\fbox{$%
\begin{array}
[c]{ccc}
{\includegraphics[
height=0.3269in,
width=0.3269in
]%
{ut05.pdf}%
}%
& {\includegraphics[
height=0.3269in,
width=0.3269in
]%
{ut05.pdf}%
}%
 \\
{\includegraphics[
height=0.3269in,
width=0.3269in
]%
{ut05.pdf}%
}%
& 
{\includegraphics[
height=0.3269in,
width=0.3269in
]%
{ut05.pdf}%
}%
\end{array}$}}
\quad\underset{K_{11}=\text{$0_1+0_1$}}{\fbox{$%
\begin{array}
[c]{ccc}
{\includegraphics[
height=0.3269in,
width=0.3269in
]%
{ut06.pdf}%
}%
& {\includegraphics[
height=0.3269in,
width=0.3269in
]%
{ut06.pdf}%
}%
 \\
{\includegraphics[
height=0.3269in,
width=0.3269in
]%
{ut06.pdf}%
}%
& 
{\includegraphics[
height=0.3269in,
width=0.3269in
]%
{ut06.pdf}%
}%
\end{array}
$}}$
$\bigskip$ \\
$\underset{K_{12}=\text{$0_1$}}{\fbox{$%
\begin{array}
[c]{ccc}
{\includegraphics[
height=0.3269in,
width=0.3269in
]%
{ut01.pdf}%
}%
& {\includegraphics[
height=0.3269in,
width=0.3269in
]%
{ut02.pdf}%
}%
 \\
{\includegraphics[
height=0.3269in,
width=0.3269in
]%
{ut03.pdf}%
}%
& 
{\includegraphics[
height=0.3269in,
width=0.3269in
]%
{ut04.pdf}%
}%
\end{array}
$}}
\quad\underset{K_{13}=\text{$0_1$}}{\fbox{$%
\begin{array}
[c]{ccc}
{\includegraphics[
height=0.3269in,
width=0.3269in
]%
{ut02.pdf}%
}%
& {\includegraphics[
height=0.3269in,
width=0.3269in
]%
{ut01.pdf}%
}%
 \\
{\includegraphics[
height=0.3269in,
width=0.3269in
]%
{ut04.pdf}%
}%
& 
{\includegraphics[
height=0.3269in,
width=0.3269in
]%
{ut03.pdf}%
}%
\end{array}
$}}
\quad\underset{K_{14}=\text{$0_1$}}{\fbox{$%
\begin{array}
[c]{ccc}
{\includegraphics[
height=0.3269in,
width=0.3269in
]%
{ut02.pdf}%
}%
& {\includegraphics[
height=0.3269in,
width=0.3269in
]%
{ut04.pdf}%
}%
 \\
{\includegraphics[
height=0.3269in,
width=0.3269in
]%
{ut03.pdf}%
}%
& 
{\includegraphics[
height=0.3269in,
width=0.3269in
]%
{ut01.pdf}%
}%
\end{array}
$}}
\quad\underset{K_{15}=\text{$0_1$}}{\fbox{$%
\begin{array}
[c]{ccc}
{\includegraphics[
height=0.3269in,
width=0.3269in
]%
{ut03.pdf}%
}%
& {\includegraphics[
height=0.3269in,
width=0.3269in
]%
{ut01.pdf}%
}%
 \\
{\includegraphics[
height=0.3269in,
width=0.3269in
]%
{ut02.pdf}%
}%
& 
{\includegraphics[
height=0.3269in,
width=0.3269in
]%
{ut04.pdf}%
}%
\end{array}$}}
\quad\underset{K_{16}=\text{$0_1$}}{\fbox{$%
\begin{array}
[c]{ccc}
{\includegraphics[
height=0.3269in,
width=0.3269in
]%
{ut03.pdf}%
}%
& {\includegraphics[
height=0.3269in,
width=0.3269in
]%
{ut01.pdf}%
}%
 \\
{\includegraphics[
height=0.3269in,
width=0.3269in
]%
{ut01.pdf}%
}%
& 
{\includegraphics[
height=0.3269in,
width=0.3269in
]%
{ut03.pdf}%
}%
\end{array}
$}}
\\ \\ \\
\underset{K_{17}=\text{$0_1$}}{\fbox{$%
\begin{array}
[c]{ccc}
{\includegraphics[
height=0.3269in,
width=0.3269in
]%
{ut02.pdf}%
}%
& {\includegraphics[
height=0.3269in,
width=0.3269in
]%
{ut04.pdf}%
}%
 \\
{\includegraphics[
height=0.3269in,
width=0.3269in
]%
{ut04.pdf}%
}%
& 
{\includegraphics[
height=0.3269in,
width=0.3269in
]%
{ut02.pdf}%
}%
\end{array}
$}}
\quad\underset{K_{18}=\text{$0_1$}}{\fbox{$%
\begin{array}
[c]{ccc}
{\includegraphics[
height=0.3269in,
width=0.3269in
]%
{ut05.pdf}%
}%
& {\includegraphics[
height=0.3269in,
width=0.3269in
]%
{ut07.pdf}%
}%
 \\
{\includegraphics[
height=0.3269in,
width=0.3269in
]%
{ut00.pdf}%
}%
& 
{\includegraphics[
height=0.3269in,
width=0.3269in
]%
{ut06.pdf}%
}%
\end{array}
$}}
\quad\underset{K_{19}=\text{$0_1$}}{\fbox{$%
\begin{array}
[c]{ccc}
{\includegraphics[
height=0.3269in,
width=0.3269in
]%
{ut00.pdf}%
}%
& {\includegraphics[
height=0.3269in,
width=0.3269in
]%
{ut06.pdf}%
}%
 \\
{\includegraphics[
height=0.3269in,
width=0.3269in
]%
{ut05.pdf}%
}%
& 
{\includegraphics[
height=0.3269in,
width=0.3269in
]%
{ut08.pdf}%
}%
\end{array}
$}}
\quad\underset{K_{20}=\text{$0_1$}}{\fbox{$%
\begin{array}
[c]{ccc}
{\includegraphics[
height=0.3269in,
width=0.3269in
]%
{ut05.pdf}%
}%
& {\includegraphics[
height=0.3269in,
width=0.3269in
]%
{ut07.pdf}%
}%
 \\
{\includegraphics[
height=0.3269in,
width=0.3269in
]%
{ut05.pdf}%
}%
& 
{\includegraphics[
height=0.3269in,
width=0.3269in
]%
{ut07.pdf}%
}%
\end{array}$}}
\quad\underset{K_{21}=\text{$0_1$}}{\fbox{$%
\begin{array}
[c]{ccc}
{\includegraphics[
height=0.3269in,
width=0.3269in
]%
{ut05.pdf}%
}%
& {\includegraphics[
height=0.3269in,
width=0.3269in
]%
{ut08.pdf}%
}%
 \\
{\includegraphics[
height=0.3269in,
width=0.3269in
]%
{ut05.pdf}%
}%
& 
{\includegraphics[
height=0.3269in,
width=0.3269in
]%
{ut08.pdf}%
}%
\end{array}
$}}
\\ \\ \\
$\bigskip$
\underset{K_{22}=\text{$0_1$}}{\fbox{$%
\begin{array}
[c]{ccc}
{\includegraphics[
height=0.3269in,
width=0.3269in
]%
{ut06.pdf}%
}%
& {\includegraphics[
height=0.3269in,
width=0.3269in
]%
{ut06.pdf}%
}%
 \\
{\includegraphics[
height=0.3269in,
width=0.3269in
]%
{ut08.pdf}%
}%
& 
{\includegraphics[
height=0.3269in,
width=0.3269in
]%
{ut08.pdf}%
}%
\end{array}
$}}
\quad\underset{K_{23}=\text{$0_1$}}{\fbox{$%
\begin{array}
[c]{ccc}
{\includegraphics[
height=0.3269in,
width=0.3269in
]%
{ut06.pdf}%
}%
& {\includegraphics[
height=0.3269in,
width=0.3269in
]%
{ut06.pdf}%
}%
 \\
{\includegraphics[
height=0.3269in,
width=0.3269in
]%
{ut07.pdf}%
}%
& 
{\includegraphics[
height=0.3269in,
width=0.3269in
]%
{ut07.pdf}%
}%
\end{array}
$}}
\quad\underset{K_{24}=\text{$0_1$}}{\fbox{$%
\begin{array}
[c]{ccc}
{\includegraphics[
height=0.3269in,
width=0.3269in
]%
{ut02.pdf}%
}%
& {\includegraphics[
height=0.3269in,
width=0.3269in
]%
{ut01.pdf}%
}%
 \\
{\includegraphics[
height=0.3269in,
width=0.3269in
]%
{ut03.pdf}%
}%
& 
{\includegraphics[
height=0.3269in,
width=0.3269in
]%
{ut04.pdf}%
}%
\end{array}
$}}
\quad\underset{K_{25}=\text{$0_1+0_1$}}{\fbox{$%
\begin{array}
[c]{ccc}
{\includegraphics[
height=0.3269in,
width=0.3269in
]%
{ut08.pdf}%
}%
& {\includegraphics[
height=0.3269in,
width=0.3269in
]%
{ut07.pdf}%
}%
 \\
{\includegraphics[
height=0.3269in,
width=0.3269in
]%
{ut07.pdf}%
}%
& 
{\includegraphics[
height=0.3269in,
width=0.3269in
]%
{ut08.pdf}%
}%
\end{array}$}}
\quad\underset{K_{26}=\text{$2_1^2$ (L2a1)}}{\fbox{$%
\begin{array}
[c]{ccc}
{\includegraphics[
height=0.3269in,
width=0.3269in
]%
{ut08.pdf}%
}%
& {\includegraphics[
height=0.3269in,
width=0.3269in
]%
{ut08.pdf}%
}%
 \\
{\includegraphics[
height=0.3269in,
width=0.3269in
]%
{ut08.pdf}%
}%
& 
{\includegraphics[
height=0.3269in,
width=0.3269in
]%
{ut08.pdf}%
}%
\end{array}
$}}
\\ \\ \\
$\bigskip$
\underset{K_{27}=\text{$2_1^2$ (L2a1)}}{\fbox{$%
\begin{array}
[c]{ccc}
{\includegraphics[
height=0.3269in,
width=0.3269in
]%
{ut07.pdf}%
}%
& {\includegraphics[
height=0.3269in,
width=0.3269in
]%
{ut07.pdf}%
}%
 \\
{\includegraphics[
height=0.3269in,
width=0.3269in
]%
{ut07.pdf}%
}%
& 
{\includegraphics[
height=0.3269in,
width=0.3269in
]%
{ut07.pdf}%
}%
\end{array}
$}}
\quad\underset{K_{28}=\text{$0_1+0_1$}}{\fbox{$%
\begin{array}
[c]{ccc}
{\includegraphics[
height=0.3269in,
width=0.3269in
]%
{ut07.pdf}%
}%
& {\includegraphics[
height=0.3269in,
width=0.3269in
]%
{ut08.pdf}%
}%
 \\
{\includegraphics[
height=0.3269in,
width=0.3269in
]%
{ut07.pdf}%
}%
& 
{\includegraphics[
height=0.3269in,
width=0.3269in
]%
{ut08.pdf}%
}%
\end{array}
$}}
\quad\underset{K_{29}=\text{$0_1+0_1$}}{\fbox{$%
\begin{array}
[c]{ccc}
{\includegraphics[
height=0.3269in,
width=0.3269in
]%
{ut08.pdf}%
}%
& {\includegraphics[
height=0.3269in,
width=0.3269in
]%
{ut07.pdf}%
}%
 \\
{\includegraphics[
height=0.3269in,
width=0.3269in
]%
{ut08.pdf}%
}%
& 
{\includegraphics[
height=0.3269in,
width=0.3269in
]%
{ut07.pdf}%
}%
\end{array}
$}}
\quad\underset{K_{30}=\text{$0_1+0_1$}}{\fbox{$%
\begin{array}
[c]{ccc}
{\includegraphics[
height=0.3269in,
width=0.3269in
]%
{ut08.pdf}%
}%
& {\includegraphics[
height=0.3269in,
width=0.3269in
]%
{ut08.pdf}%
}%
 \\
{\includegraphics[
height=0.3269in,
width=0.3269in
]%
{ut07.pdf}%
}%
& 
{\includegraphics[
height=0.3269in,
width=0.3269in
]%
{ut07.pdf}%
}%
\end{array}$}}
\quad\underset{K_{31}=\text{$0_1+0_1$}}{\fbox{$%
\begin{array}
[c]{ccc}
{\includegraphics[
height=0.3269in,
width=0.3269in
]%
{ut07.pdf}%
}%
& {\includegraphics[
height=0.3269in,
width=0.3269in
]%
{ut07.pdf}%
}%
 \\
{\includegraphics[
height=0.3269in,
width=0.3269in
]%
{ut08.pdf}%
}%
& 
{\includegraphics[
height=0.3269in,
width=0.3269in
]%
{ut08.pdf}%
}%
\end{array}
$}}
\\ \\ \\
$\bigskip$
\underset{K_{32}=\text{$0_1$}}{\fbox{$%
\begin{array}
[c]{ccc}
{\includegraphics[
height=0.3269in,
width=0.3269in
]%
{ut08.pdf}%
}%
& {\includegraphics[
height=0.3269in,
width=0.3269in
]%
{ut07.pdf}%
}%
 \\
{\includegraphics[
height=0.3269in,
width=0.3269in
]%
{ut08.pdf}%
}%
& 
{\includegraphics[
height=0.3269in,
width=0.3269in
]%
{ut08.pdf}%
}%
\end{array}
$}}
\quad\underset{K_{33}=\text{$2_1^2$ (L2a1)}}{\fbox{$%
\begin{array}
[c]{ccc}
{\includegraphics[
height=0.3269in,
width=0.3269in
]%
{ut05.pdf}%
}%
& {\includegraphics[
height=0.3269in,
width=0.3269in
]%
{ut10.pdf}%
}%
 \\
{\includegraphics[
height=0.3269in,
width=0.3269in
]%
{ut00.pdf}%
}%
& 
{\includegraphics[
height=0.3269in,
width=0.3269in
]%
{ut06.pdf}%
}%
\end{array}
$}}
\quad\underset{K_{34}=\text{$0_1+0_1$}}{\fbox{$%
\begin{array}
[c]{ccc}
{\includegraphics[
height=0.3269in,
width=0.3269in
]%
{ut05.pdf}%
}%
& {\includegraphics[
height=0.3269in,
width=0.3269in
]%
{ut09.pdf}%
}%
 \\
{\includegraphics[
height=0.3269in,
width=0.3269in
]%
{ut00.pdf}%
}%
& 
{\includegraphics[
height=0.3269in,
width=0.3269in
]%
{ut06.pdf}%
}%
\end{array}
$}}
\quad\underset{K_{35}=\text{$0_1$}}{\fbox{$%
\begin{array}
[c]{ccc}
{\includegraphics[
height=0.3269in,
width=0.3269in
]%
{ut08.pdf}%
}%
& {\includegraphics[
height=0.3269in,
width=0.3269in
]%
{ut10.pdf}%
}%
 \\
{\includegraphics[
height=0.3269in,
width=0.3269in
]%
{ut07.pdf}%
}%
& 
{\includegraphics[
height=0.3269in,
width=0.3269in
]%
{ut08.pdf}%
}%
\end{array}$}}
\quad\underset{K_{36}=\text{$0_1$}}{\fbox{$%
\begin{array}
[c]{ccc}
{\includegraphics[
height=0.3269in,
width=0.3269in
]%
{ut08.pdf}%
}%
& {\includegraphics[
height=0.3269in,
width=0.3269in
]%
{ut09.pdf}%
}%
 \\
{\includegraphics[
height=0.3269in,
width=0.3269in
]%
{ut07.pdf}%
}%
& 
{\includegraphics[
height=0.3269in,
width=0.3269in
]%
{ut08.pdf}%
}%
\end{array}
$}}
\\ \\ \\
$\bigskip$
\underset{K_{37}=\text{$0_1$}}{\fbox{$%
\begin{array}
[c]{ccc}
{\includegraphics[
height=0.3269in,
width=0.3269in
]%
{ut08.pdf}%
}%
& {\includegraphics[
height=0.3269in,
width=0.3269in
]%
{ut07.pdf}%
}%
 \\
{\includegraphics[
height=0.3269in,
width=0.3269in
]%
{ut07.pdf}%
}%
& 
{\includegraphics[
height=0.3269in,
width=0.3269in
]%
{ut09.pdf}%
}%
\end{array}
$}}
\quad\underset{K_{38}=\text{$0_1+0_1$}}{\fbox{$%
\begin{array}
[c]{ccc}
{\includegraphics[
height=0.3269in,
width=0.3269in
]%
{ut08.pdf}%
}%
& {\includegraphics[
height=0.3269in,
width=0.3269in
]%
{ut07.pdf}%
}%
 \\
{\includegraphics[
height=0.3269in,
width=0.3269in
]%
{ut07.pdf}%
}%
& 
{\includegraphics[
height=0.3269in,
width=0.3269in
]%
{ut10.pdf}%
}%
\end{array}
$}}
\quad\underset{K_{39}=\text{$2_1^2$ (L2a1)}}{\fbox{$%
\begin{array}
[c]{ccc}
{\includegraphics[
height=0.3269in,
width=0.3269in
]%
{ut05.pdf}%
}%
& {\includegraphics[
height=0.3269in,
width=0.3269in
]%
{ut10.pdf}%
}%
 \\
{\includegraphics[
height=0.3269in,
width=0.3269in
]%
{ut05.pdf}%
}%
& 
{\includegraphics[
height=0.3269in,
width=0.3269in
]%
{ut08.pdf}%
}%
\end{array}
$}}
\quad\underset{K_{40}=\text{$0_1+0_1$}}{\fbox{$%
\begin{array}
[c]{ccc}
{\includegraphics[
height=0.3269in,
width=0.3269in
]%
{ut05.pdf}%
}%
& {\includegraphics[
height=0.3269in,
width=0.3269in
]%
{ut09.pdf}%
}%
 \\
{\includegraphics[
height=0.3269in,
width=0.3269in
]%
{ut05.pdf}%
}%
& 
{\includegraphics[
height=0.3269in,
width=0.3269in
]%
{ut08.pdf}%
}%
\end{array}$}}
\quad\underset{K_{41}=\text{$2_1^2$ (L2a1)}}{\fbox{$%
\begin{array}
[c]{ccc}
{\includegraphics[
height=0.3269in,
width=0.3269in
]%
{ut06.pdf}%
}%
& {\includegraphics[
height=0.3269in,
width=0.3269in
]%
{ut06.pdf}%
}%
 \\
{\includegraphics[
height=0.3269in,
width=0.3269in
]%
{ut08.pdf}%
}%
& 
{\includegraphics[
height=0.3269in,
width=0.3269in
]%
{ut10.pdf}%
}%
\end{array}
$}}
\\ \\ \\
$\bigskip$
\underset{K_{42}=\text{$0_1+0_1$}}{\fbox{$%
\begin{array}
[c]{ccc}
{\includegraphics[
height=0.3269in,
width=0.3269in
]%
{ut06.pdf}%
}%
& {\includegraphics[
height=0.3269in,
width=0.3269in
]%
{ut06.pdf}%
}%
 \\
{\includegraphics[
height=0.3269in,
width=0.3269in
]%
{ut08.pdf}%
}%
& 
{\includegraphics[
height=0.3269in,
width=0.3269in
]%
{ut09.pdf}%
}%
\end{array}
$}}
\quad\underset{K_{43}=\text{$2_1^2$ (L2a1)}}{\fbox{$%
\begin{array}
[c]{ccc}
{\includegraphics[
height=0.3269in,
width=0.3269in
]%
{ut05.pdf}%
}%
& {\includegraphics[
height=0.3269in,
width=0.3269in
]%
{ut10.pdf}%
}%
 \\
{\includegraphics[
height=0.3269in,
width=0.3269in
]%
{ut05.pdf}%
}%
& 
{\includegraphics[
height=0.3269in,
width=0.3269in
]%
{ut07.pdf}%
}%
\end{array}
$}}
\quad\underset{K_{44}=\text{$0_1+0_1$}}{\fbox{$%
\begin{array}
[c]{ccc}
{\includegraphics[
height=0.3269in,
width=0.3269in
]%
{ut05.pdf}%
}%
& {\includegraphics[
height=0.3269in,
width=0.3269in
]%
{ut09.pdf}%
}%
 \\
{\includegraphics[
height=0.3269in,
width=0.3269in
]%
{ut05.pdf}%
}%
& 
{\includegraphics[
height=0.3269in,
width=0.3269in
]%
{ut07.pdf}%
}%
\end{array}
$}}
\quad\underset{K_{45}=\text{$2_1^2$ (L2a1)}}{\fbox{$%
\begin{array}
[c]{ccc}
{\includegraphics[
height=0.3269in,
width=0.3269in
]%
{ut06.pdf}%
}%
& {\includegraphics[
height=0.3269in,
width=0.3269in
]%
{ut06.pdf}%
}%
 \\
{\includegraphics[
height=0.3269in,
width=0.3269in
]%
{ut07.pdf}%
}%
& 
{\includegraphics[
height=0.3269in,
width=0.3269in
]%
{ut10.pdf}%
}%
\end{array}$}}
\quad\underset{K_{46}=\text{$0_1+0_1$}}{\fbox{$%
\begin{array}
[c]{ccc}
{\includegraphics[
height=0.3269in,
width=0.3269in
]%
{ut06.pdf}%
}%
& {\includegraphics[
height=0.3269in,
width=0.3269in
]%
{ut06.pdf}%
}%
 \\
{\includegraphics[
height=0.3269in,
width=0.3269in
]%
{ut07.pdf}%
}%
& 
{\includegraphics[
height=0.3269in,
width=0.3269in
]%
{ut09.pdf}%
}%
\end{array}
$}}
\\ \\ \\
$\bigskip$
\underset{K_{47}=\text{$0_1$}}{\fbox{$%
\begin{array}
[c]{ccc}
{\includegraphics[
height=0.3269in,
width=0.3269in
]%
{ut07.pdf}%
}%
& {\includegraphics[
height=0.3269in,
width=0.3269in
]%
{ut10.pdf}%
}%
 \\
{\includegraphics[
height=0.3269in,
width=0.3269in
]%
{ut08.pdf}%
}%
& 
{\includegraphics[
height=0.3269in,
width=0.3269in
]%
{ut08.pdf}%
}%
\end{array}
$}}
\quad\underset{K_{48}=\text{$0_1$}}{\fbox{$%
\begin{array}
[c]{ccc}
{\includegraphics[
height=0.3269in,
width=0.3269in
]%
{ut07.pdf}%
}%
& {\includegraphics[
height=0.3269in,
width=0.3269in
]%
{ut09.pdf}%
}%
 \\
{\includegraphics[
height=0.3269in,
width=0.3269in
]%
{ut08.pdf}%
}%
& 
{\includegraphics[
height=0.3269in,
width=0.3269in
]%
{ut08.pdf}%
}%
\end{array}
$}}
\quad\underset{K_{49}=\text{$3_1$}}{\fbox{$%
\begin{array}
[c]{ccc}
{\includegraphics[
height=0.3269in,
width=0.3269in
]%
{ut07.pdf}%
}%
& {\includegraphics[
height=0.3269in,
width=0.3269in
]%
{ut10.pdf}%
}%
 \\
{\includegraphics[
height=0.3269in,
width=0.3269in
]%
{ut07.pdf}%
}%
& 
{\includegraphics[
height=0.3269in,
width=0.3269in
]%
{ut07.pdf}%
}%
\end{array}
$}}
\quad\underset{K_{50}=\text{$0_1$}}{\fbox{$%
\begin{array}
[c]{ccc}
{\includegraphics[
height=0.3269in,
width=0.3269in
]%
{ut07.pdf}%
}%
& {\includegraphics[
height=0.3269in,
width=0.3269in
]%
{ut09.pdf}%
}%
 \\
{\includegraphics[
height=0.3269in,
width=0.3269in
]%
{ut07.pdf}%
}%
& 
{\includegraphics[
height=0.3269in,
width=0.3269in
]%
{ut07.pdf}%
}%
\end{array}$}}
\quad\underset{K_{51}=\text{$0_1$}}{\fbox{$%
\begin{array}
[c]{ccc}
{\includegraphics[
height=0.3269in,
width=0.3269in
]%
{ut07.pdf}%
}%
& {\includegraphics[
height=0.3269in,
width=0.3269in
]%
{ut10.pdf}%
}%
 \\
{\includegraphics[
height=0.3269in,
width=0.3269in
]%
{ut07.pdf}%
}%
& 
{\includegraphics[
height=0.3269in,
width=0.3269in
]%
{ut08.pdf}%
}%
\end{array}
$}}
\\ \\ \\
$\bigskip$
\underset{K_{52}=\text{$0_1$}}{\fbox{$%
\begin{array}
[c]{ccc}
{\includegraphics[
height=0.3269in,
width=0.3269in
]%
{ut07.pdf}%
}%
& {\includegraphics[
height=0.3269in,
width=0.3269in
]%
{ut09.pdf}%
}%
 \\
{\includegraphics[
height=0.3269in,
width=0.3269in
]%
{ut07.pdf}%
}%
& 
{\includegraphics[
height=0.3269in,
width=0.3269in
]%
{ut08.pdf}%
}%
\end{array}
$}}
\quad\underset{K_{53}=\text{$3_1$}}{\fbox{$%
\begin{array}
[c]{ccc}
{\includegraphics[
height=0.3269in,
width=0.3269in
]%
{ut08.pdf}%
}%
& {\includegraphics[
height=0.3269in,
width=0.3269in
]%
{ut10.pdf}%
}%
 \\
{\includegraphics[
height=0.3269in,
width=0.3269in
]%
{ut08.pdf}%
}%
& 
{\includegraphics[
height=0.3269in,
width=0.3269in
]%
{ut08.pdf}%
}%
\end{array}
$}}
\quad\underset{K_{54}=\text{$0_1$}}{\fbox{$%
\begin{array}
[c]{ccc}
{\includegraphics[
height=0.3269in,
width=0.3269in
]%
{ut08.pdf}%
}%
& {\includegraphics[
height=0.3269in,
width=0.3269in
]%
{ut09.pdf}%
}%
 \\
{\includegraphics[
height=0.3269in,
width=0.3269in
]%
{ut08.pdf}%
}%
& 
{\includegraphics[
height=0.3269in,
width=0.3269in
]%
{ut08.pdf}%
}%
\end{array}
$}}
\quad\underset{K_{55}=\text{$0_1$}}{\fbox{$%
\begin{array}
[c]{ccc}
{\includegraphics[
height=0.3269in,
width=0.3269in
]%
{ut08.pdf}%
}%
& {\includegraphics[
height=0.3269in,
width=0.3269in
]%
{ut10.pdf}%
}%
 \\
{\includegraphics[
height=0.3269in,
width=0.3269in
]%
{ut08.pdf}%
}%
& 
{\includegraphics[
height=0.3269in,
width=0.3269in
]%
{ut07.pdf}%
}%
\end{array}$}}
\quad\underset{K_{56}=\text{$0_1$}}{\fbox{$%
\begin{array}
[c]{ccc}
{\includegraphics[
height=0.3269in,
width=0.3269in
]%
{ut08.pdf}%
}%
& {\includegraphics[
height=0.3269in,
width=0.3269in
]%
{ut09.pdf}%
}%
 \\
{\includegraphics[
height=0.3269in,
width=0.3269in
]%
{ut08.pdf}%
}%
& 
{\includegraphics[
height=0.3269in,
width=0.3269in
]%
{ut07.pdf}%
}%
\end{array}
$}}
\\ \\ \\
$\bigskip$
\underset{K_{57}=\text{$0_1$}}{\fbox{$%
\begin{array}
[c]{ccc}
{\includegraphics[
height=0.3269in,
width=0.3269in
]%
{ut07.pdf}%
}%
& {\includegraphics[
height=0.3269in,
width=0.3269in
]%
{ut07.pdf}%
}%
 \\
{\includegraphics[
height=0.3269in,
width=0.3269in
]%
{ut08.pdf}%
}%
& 
{\includegraphics[
height=0.3269in,
width=0.3269in
]%
{ut10.pdf}%
}%
\end{array}
$}}
\quad\underset{K_{58}=\text{$0_1$}}{\fbox{$%
\begin{array}
[c]{ccc}
{\includegraphics[
height=0.3269in,
width=0.3269in
]%
{ut07.pdf}%
}%
& {\includegraphics[
height=0.3269in,
width=0.3269in
]%
{ut07.pdf}%
}%
 \\
{\includegraphics[
height=0.3269in,
width=0.3269in
]%
{ut08.pdf}%
}%
& 
{\includegraphics[
height=0.3269in,
width=0.3269in
]%
{ut09.pdf}%
}%
\end{array}
$}}
\quad\underset{K_{59}=\text{$0_1+0_1+0_1$}}{\fbox{$%
\begin{array}
[c]{ccc}
{\includegraphics[
height=0.3269in,
width=0.3269in
]%
{ut05.pdf}%
}%
& {\includegraphics[
height=0.3269in,
width=0.3269in
]%
{ut09.pdf}%
}%
 \\
{\includegraphics[
height=0.3269in,
width=0.3269in
]%
{ut05.pdf}%
}%
& 
{\includegraphics[
height=0.3269in,
width=0.3269in
]%
{ut09.pdf}%
}%
\end{array}
$}}
\quad\underset{K_{60}=\text{$2_1^2+0_1$}}{\fbox{$%
\begin{array}
[c]{ccc}
{\includegraphics[
height=0.3269in,
width=0.3269in
]%
{ut05.pdf}%
}%
& {\includegraphics[
height=0.3269in,
width=0.3269in
]%
{ut09.pdf}%
}%
 \\
{\includegraphics[
height=0.3269in,
width=0.3269in
]%
{ut05.pdf}%
}%
& 
{\includegraphics[
height=0.3269in,
width=0.3269in
]%
{ut10.pdf}%
}%
\end{array}
$}}
\quad\underset{K_{61}=\text{$2_1^2\#2_1^2$}}{\fbox{$%
\begin{array}
[c]{ccc}
{\includegraphics[
height=0.3269in,
width=0.3269in
]%
{ut05.pdf}%
}%
& {\includegraphics[
height=0.3269in,
width=0.3269in
]%
{ut10.pdf}%
}%
 \\
{\includegraphics[
height=0.3269in,
width=0.3269in
]%
{ut05.pdf}%
}%
& 
{\includegraphics[
height=0.3269in,
width=0.3269in
]%
{ut10.pdf}%
}%
\end{array}
$}}
\\ \\ \\
$\bigskip$
\underset{K_{62}=\text{$0_1+0_1+0_1$}}{\fbox{$%
\begin{array}
[c]{ccc}
{\includegraphics[
height=0.3269in,
width=0.3269in
]%
{ut06.pdf}%
}%
& {\includegraphics[
height=0.3269in,
width=0.3269in
]%
{ut06.pdf}%
}%
 \\
{\includegraphics[
height=0.3269in,
width=0.3269in
]%
{ut09.pdf}%
}%
& 
{\includegraphics[
height=0.3269in,
width=0.3269in
]%
{ut09.pdf}%
}%
\end{array}
$}}
\quad\underset{K_{63}=\text{$2_1^2+0_1$}}{\fbox{$%
\begin{array}
[c]{ccc}
{\includegraphics[
height=0.3269in,
width=0.3269in
]%
{ut06.pdf}%
}%
& {\includegraphics[
height=0.3269in,
width=0.3269in
]%
{ut06.pdf}%
}%
 \\
{\includegraphics[
height=0.3269in,
width=0.3269in
]%
{ut09.pdf}%
}%
& 
{\includegraphics[
height=0.3269in,
width=0.3269in
]%
{ut10.pdf}%
}%
\end{array}$}}
\quad\underset{K_{64}=\text{$2_1^2\#2_1^2$ }}{\fbox{$%
\begin{array}
[c]{ccc}
{\includegraphics[
height=0.3269in,
width=0.3269in
]%
{ut06.pdf}%
}%
& {\includegraphics[
height=0.3269in,
width=0.3269in
]%
{ut06.pdf}%
}%
 \\
{\includegraphics[
height=0.3269in,
width=0.3269in
]%
{ut10.pdf}%
}%
& 
{\includegraphics[
height=0.3269in,
width=0.3269in
]%
{ut10.pdf}%
}%
\end{array}
$}}
\quad\underset{K_{65}=\text{$0_1+0_1$}}{\fbox{$%
\begin{array}
[c]{ccc}
{\includegraphics[
height=0.3269in,
width=0.3269in
]%
{ut08.pdf}%
}%
& {\includegraphics[
height=0.3269in,
width=0.3269in
]%
{ut10.pdf}%
}%
 \\
{\includegraphics[
height=0.3269in,
width=0.3269in
]%
{ut07.pdf}%
}%
& 
{\includegraphics[
height=0.3269in,
width=0.3269in
]%
{ut10.pdf}%
}%
\end{array}
$}}
\quad\underset{K_{66}=\text{$2_1^2$ (L2a1)}}{\fbox{$%
\begin{array}
[c]{ccc}
{\includegraphics[
height=0.3269in,
width=0.3269in
]%
{ut08.pdf}%
}%
& {\includegraphics[
height=0.3269in,
width=0.3269in
]%
{ut09.pdf}%
}%
 \\
{\includegraphics[
height=0.3269in,
width=0.3269in
]%
{ut07.pdf}%
}%
& 
{\includegraphics[
height=0.3269in,
width=0.3269in
]%
{ut10.pdf}%
}%
\end{array}
$}}
\\ \\ \\
$\bigskip$
\underset{K_{67}=\text{$0_1+0_1$}}{\fbox{$%
\begin{array}
[c]{ccc}
{\includegraphics[
height=0.3269in,
width=0.3269in
]%
{ut08.pdf}%
}%
& {\includegraphics[
height=0.3269in,
width=0.3269in
]%
{ut09.pdf}%
}%
 \\
{\includegraphics[
height=0.3269in,
width=0.3269in
]%
{ut07.pdf}%
}%
& 
{\includegraphics[
height=0.3269in,
width=0.3269in
]%
{ut09.pdf}%
}%
\end{array}
$}}
\quad\underset{K_{68}=\text{$0_1+0_1$}}{\fbox{$%
\begin{array}
[c]{ccc}
{\includegraphics[
height=0.3269in,
width=0.3269in
]%
{ut08.pdf}%
}%
& {\includegraphics[
height=0.3269in,
width=0.3269in
]%
{ut07.pdf}%
}%
 \\
{\includegraphics[
height=0.3269in,
width=0.3269in
]%
{ut09.pdf}%
}%
& 
{\includegraphics[
height=0.3269in,
width=0.3269in
]%
{ut09.pdf}%
}%
\end{array}$}}
\quad\underset{K_{69}=\text{$2_1^2$ (L2a1)}}{\fbox{$%
\begin{array}
[c]{ccc}
{\includegraphics[
height=0.3269in,
width=0.3269in
]%
{ut08.pdf}%
}%
& {\includegraphics[
height=0.3269in,
width=0.3269in
]%
{ut07.pdf}%
}%
 \\
{\includegraphics[
height=0.3269in,
width=0.3269in
]%
{ut09.pdf}%
}%
& 
{\includegraphics[
height=0.3269in,
width=0.3269in
]%
{ut10.pdf}%
}%
\end{array}
$}}
\quad\underset{K_{70}=\text{$0_1+0_1$}}{\fbox{$%
\begin{array}
[c]{ccc}
{\includegraphics[
height=0.3269in,
width=0.3269in
]%
{ut08.pdf}%
}%
& {\includegraphics[
height=0.3269in,
width=0.3269in
]%
{ut07.pdf}%
}%
 \\
{\includegraphics[
height=0.3269in,
width=0.3269in
]%
{ut10.pdf}%
}%
& 
{\includegraphics[
height=0.3269in,
width=0.3269in
]%
{ut10.pdf}%
}%
\end{array}
$}}
\quad\underset{K_{71}=\text{$2_1^2$ (L2a1)}}{\fbox{$%
\begin{array}
[c]{ccc}
{\includegraphics[
height=0.3269in,
width=0.3269in
]%
{ut07.pdf}%
}%
& {\includegraphics[
height=0.3269in,
width=0.3269in
]%
{ut10.pdf}%
}%
 \\
{\includegraphics[
height=0.3269in,
width=0.3269in
]%
{ut07.pdf}%
}%
& 
{\includegraphics[
height=0.3269in,
width=0.3269in
]%
{ut09.pdf}%
}%
\end{array}
$}}
\\ \\ \\
\underset{K_{72}=\text{$0_1+0_1$}}{\fbox{$%
\begin{array}
[c]{ccc}
{\includegraphics[
height=0.3269in,
width=0.3269in
]%
{ut07.pdf}%
}%
& {\includegraphics[
height=0.3269in,
width=0.3269in
]%
{ut09.pdf}%
}%
 \\
{\includegraphics[
height=0.3269in,
width=0.3269in
]%
{ut07.pdf}%
}%
& 
{\includegraphics[
height=0.3269in,
width=0.3269in
]%
{ut09.pdf}%
}%
\end{array}
$}}
\quad\underset{K_{73}=\text{$4_1^2$ (L4a1)}}{\fbox{$%
\begin{array}
[c]{ccc}
{\includegraphics[
height=0.3269in,
width=0.3269in
]%
{ut07.pdf}%
}%
& {\includegraphics[
height=0.3269in,
width=0.3269in
]%
{ut10.pdf}%
}%
 \\
{\includegraphics[
height=0.3269in,
width=0.3269in
]%
{ut07.pdf}%
}%
& 
{\includegraphics[
height=0.3269in,
width=0.3269in
]%
{ut10.pdf}%
}%
\end{array}$}}
\quad\underset{K_{74}=\text{$0_1+0_1$}}{\fbox{$%
\begin{array}
[c]{ccc}
{\includegraphics[
height=0.3269in,
width=0.3269in
]%
{ut08.pdf}%
}%
& {\includegraphics[
height=0.3269in,
width=0.3269in
]%
{ut08.pdf}%
}%
 \\
{\includegraphics[
height=0.3269in,
width=0.3269in
]%
{ut09.pdf}%
}%
& 
{\includegraphics[
height=0.3269in,
width=0.3269in
]%
{ut09.pdf}%
}%
\end{array}
$}}
\quad\underset{K_{75}=\text{$4_1^2$ (L4a1)}}{\fbox{$%
\begin{array}
[c]{ccc}
{\includegraphics[
height=0.3269in,
width=0.3269in
]%
{ut08.pdf}%
}%
& {\includegraphics[
height=0.3269in,
width=0.3269in
]%
{ut08.pdf}%
}%
 \\
{\includegraphics[
height=0.3269in,
width=0.3269in
]%
{ut10.pdf}%
}%
& 
{\includegraphics[
height=0.3269in,
width=0.3269in
]%
{ut10.pdf}%
}%
\end{array}
$}}
\quad\underset{K_{76}=\text{$2_1^2$ (L2a1)}}{\fbox{$%
\begin{array}
[c]{ccc}
{\includegraphics[
height=0.3269in,
width=0.3269in
]%
{ut08.pdf}%
}%
& {\includegraphics[
height=0.3269in,
width=0.3269in
]%
{ut08.pdf}%
}%
 \\
{\includegraphics[
height=0.3269in,
width=0.3269in
]%
{ut09.pdf}%
}%
& 
{\includegraphics[
height=0.3269in,
width=0.3269in
]%
{ut10.pdf}%
}%
\end{array}
$}}
\\ \\ \\
$\bigskip$
\underset{K_{77}=\text{$2_1^2$ (L2a1)}}{\fbox{$%
\begin{array}
[c]{ccc}
{\includegraphics[
height=0.3269in,
width=0.3269in
]%
{ut07.pdf}%
}%
& {\includegraphics[
height=0.3269in,
width=0.3269in
]%
{ut10.pdf}%
}%
 \\
{\includegraphics[
height=0.3269in,
width=0.3269in
]%
{ut08.pdf}%
}%
& 
{\includegraphics[
height=0.3269in,
width=0.3269in
]%
{ut09.pdf}%
}%
\end{array}
$}}
\quad\underset{K_{78}=\text{$0_1+0_1$}}{\fbox{$%
\begin{array}
[c]{ccc}
{\includegraphics[
height=0.3269in,
width=0.3269in
]%
{ut07.pdf}%
}%
& {\includegraphics[
height=0.3269in,
width=0.3269in
]%
{ut10.pdf}%
}%
 \\
{\includegraphics[
height=0.3269in,
width=0.3269in
]%
{ut08.pdf}%
}%
& 
{\includegraphics[
height=0.3269in,
width=0.3269in
]%
{ut10.pdf}%
}%
\end{array}$}}
\quad\underset{K_{79}=\text{$0_1+0_1$}}{\fbox{$%
\begin{array}
[c]{ccc}
{\includegraphics[
height=0.3269in,
width=0.3269in
]%
{ut07.pdf}%
}%
& {\includegraphics[
height=0.3269in,
width=0.3269in
]%
{ut09.pdf}%
}%
 \\
{\includegraphics[
height=0.3269in,
width=0.3269in
]%
{ut08.pdf}%
}%
& 
{\includegraphics[
height=0.3269in,
width=0.3269in
]%
{ut09.pdf}%
}%
\end{array}
$}}
\quad\underset{K_{80}=\text{$2_1^2$ (L2a1)}}{\fbox{$%
\begin{array}
[c]{ccc}
{\includegraphics[
height=0.3269in,
width=0.3269in
]%
{ut07.pdf}%
}%
& {\includegraphics[
height=0.3269in,
width=0.3269in
]%
{ut08.pdf}%
}%
 \\
{\includegraphics[
height=0.3269in,
width=0.3269in
]%
{ut09.pdf}%
}%
& 
{\includegraphics[
height=0.3269in,
width=0.3269in
]%
{ut10.pdf}%
}%
\end{array}
$}}
\quad\underset{K_{81}=\text{$0_1+0_1$}}{\fbox{$%
\begin{array}
[c]{ccc}
{\includegraphics[
height=0.3269in,
width=0.3269in
]%
{ut07.pdf}%
}%
& {\includegraphics[
height=0.3269in,
width=0.3269in
]%
{ut08.pdf}%
}%
 \\
{\includegraphics[
height=0.3269in,
width=0.3269in
]%
{ut09.pdf}%
}%
& 
{\includegraphics[
height=0.3269in,
width=0.3269in
]%
{ut09.pdf}%
}%
\end{array}
$}}
\\ \\ \\
$\bigskip$
\underset{K_{82}=\text{$0_1+0_1$}}{\fbox{$%
\begin{array}
[c]{ccc}
{\includegraphics[
height=0.3269in,
width=0.3269in
]%
{ut07.pdf}%
}%
& {\includegraphics[
height=0.3269in,
width=0.3269in
]%
{ut08.pdf}%
}%
 \\
{\includegraphics[
height=0.3269in,
width=0.3269in
]%
{ut10.pdf}%
}%
& 
{\includegraphics[
height=0.3269in,
width=0.3269in
]%
{ut10.pdf}%
}%
\end{array}
$}}
\quad\underset{K_{83}=\text{$4_1^2$ (L4a1)}}{\fbox{$%
\begin{array}
[c]{ccc}
{\includegraphics[
height=0.3269in,
width=0.3269in
]%
{ut10.pdf}%
}%
& {\includegraphics[
height=0.3269in,
width=0.3269in
]%
{ut07.pdf}%
}%
 \\
{\includegraphics[
height=0.3269in,
width=0.3269in
]%
{ut07.pdf}%
}%
& 
{\includegraphics[
height=0.3269in,
width=0.3269in
]%
{ut10.pdf}%
}%
\end{array}$}}
\quad\underset{K_{84}=\text{$0_1+0_1$}}{\fbox{$%
\begin{array}
[c]{ccc}
{\includegraphics[
height=0.3269in,
width=0.3269in
]%
{ut10.pdf}%
}%
& {\includegraphics[
height=0.3269in,
width=0.3269in
]%
{ut07.pdf}%
}%
 \\
{\includegraphics[
height=0.3269in,
width=0.3269in
]%
{ut08.pdf}%
}%
& 
{\includegraphics[
height=0.3269in,
width=0.3269in
]%
{ut10.pdf}%
}%
\end{array}
$}}
\quad\underset{K_{85}=\text{$4_1^2$ (L4a1)}}{\fbox{$%
\begin{array}
[c]{ccc}
{\includegraphics[
height=0.3269in,
width=0.3269in
]%
{ut10.pdf}%
}%
& {\includegraphics[
height=0.3269in,
width=0.3269in
]%
{ut08.pdf}%
}%
 \\
{\includegraphics[
height=0.3269in,
width=0.3269in
]%
{ut08.pdf}%
}%
& 
{\includegraphics[
height=0.3269in,
width=0.3269in
]%
{ut10.pdf}%
}%
\end{array}
$}}
\quad\underset{K_{86}=\text{$2_1^2$ (L2a1)}}{\fbox{$%
\begin{array}
[c]{ccc}
{\includegraphics[
height=0.3269in,
width=0.3269in
]%
{ut09.pdf}%
}%
& {\includegraphics[
height=0.3269in,
width=0.3269in
]%
{ut07.pdf}%
}%
 \\
{\includegraphics[
height=0.3269in,
width=0.3269in
]%
{ut07.pdf}%
}%
& 
{\includegraphics[
height=0.3269in,
width=0.3269in
]%
{ut10.pdf}%
}%
\end{array}
$}}
\\ \\ \\
$\bigskip$
\underset{K_{87}=\text{$2_1^2$ (L2a1)}}{\fbox{$%
\begin{array}
[c]{ccc}
{\includegraphics[
height=0.3269in,
width=0.3269in
]%
{ut09.pdf}%
}%
& {\includegraphics[
height=0.3269in,
width=0.3269in
]%
{ut07.pdf}%
}%
 \\
{\includegraphics[
height=0.3269in,
width=0.3269in
]%
{ut08.pdf}%
}%
& 
{\includegraphics[
height=0.3269in,
width=0.3269in
]%
{ut10.pdf}%
}%
\end{array}
$}}
\quad\underset{K_{88}=\text{$4_1^2$ (L4a1)}}{\fbox{$%
\begin{array}
[c]{ccc}
{\includegraphics[
height=0.3269in,
width=0.3269in
]%
{ut09.pdf}%
}%
& {\includegraphics[
height=0.3269in,
width=0.3269in
]%
{ut08.pdf}%
}%
 \\
{\includegraphics[
height=0.3269in,
width=0.3269in
]%
{ut08.pdf}%
}%
& 
{\includegraphics[
height=0.3269in,
width=0.3269in
]%
{ut10.pdf}%
}%
\end{array}$}}
\quad\underset{K_{89}=\text{$2_1^2+0_1$}}{\fbox{$%
\begin{array}
[c]{ccc}
{\includegraphics[
height=0.3269in,
width=0.3269in
]%
{ut09.pdf}%
}%
& {\includegraphics[
height=0.3269in,
width=0.3269in
]%
{ut09.pdf}%
}%
 \\
{\includegraphics[
height=0.3269in,
width=0.3269in
]%
{ut07.pdf}%
}%
& 
{\includegraphics[
height=0.3269in,
width=0.3269in
]%
{ut10.pdf}%
}%
\end{array}$}}
\quad\underset{K_{90}=\text{$2_1^2+0_1$}}{\fbox{$%
\begin{array}
[c]{ccc}
{\includegraphics[
height=0.3269in,
width=0.3269in
]%
{ut09.pdf}%
}%
& {\includegraphics[
height=0.3269in,
width=0.3269in
]%
{ut09.pdf}%
}%
 \\
{\includegraphics[
height=0.3269in,
width=0.3269in
]%
{ut08.pdf}%
}%
& 
{\includegraphics[
height=0.3269in,
width=0.3269in
]%
{ut10.pdf}%
}%
\end{array}
$}}
\quad\underset{K_{91}=\text{$2_1^2\#2_1^2$}}{\fbox{$%
\begin{array}
[c]{ccc}
{\includegraphics[
height=0.3269in,
width=0.3269in
]%
{ut09.pdf}%
}%
& {\includegraphics[
height=0.3269in,
width=0.3269in
]%
{ut10.pdf}%
}%
 \\
{\includegraphics[
height=0.3269in,
width=0.3269in
]%
{ut07.pdf}%
}%
& 
{\includegraphics[
height=0.3269in,
width=0.3269in
]%
{ut10.pdf}%
}%
\end{array}
$}}
\\ \\ \\
$\bigskip$
\underset{K_{92}=\text{$2_1^2\#2_1^2$}}{\fbox{$%
\begin{array}
[c]{ccc}
{\includegraphics[
height=0.3269in,
width=0.3269in
]%
{ut09.pdf}%
}%
& {\includegraphics[
height=0.3269in,
width=0.3269in
]%
{ut10.pdf}%
}%
 \\
{\includegraphics[
height=0.3269in,
width=0.3269in
]%
{ut08.pdf}%
}%
& 
{\includegraphics[
height=0.3269in,
width=0.3269in
]%
{ut10.pdf}%
}%
\end{array}
$}}
\quad\underset{K_{93}=\text{$6_1^3$ (L6a5)}}{\fbox{$%
\begin{array}
[c]{ccc}
{\includegraphics[
height=0.3269in,
width=0.3269in
]%
{ut10.pdf}%
}%
& {\includegraphics[
height=0.3269in,
width=0.3269in
]%
{ut10.pdf}%
}%
 \\
{\includegraphics[
height=0.3269in,
width=0.3269in
]%
{ut07.pdf}%
}%
& 
{\includegraphics[
height=0.3269in,
width=0.3269in
]%
{ut10.pdf}%
}%
\end{array}$}}
\quad\underset{K_{94}=\text{$6_1^3$ (L6a5)}}{\fbox{$%
\begin{array}
[c]{ccc}
{\includegraphics[
height=0.3269in,
width=0.3269in
]%
{ut10.pdf}%
}%
& {\includegraphics[
height=0.3269in,
width=0.3269in
]%
{ut10.pdf}%
}%
 \\
{\includegraphics[
height=0.3269in,
width=0.3269in
]%
{ut08.pdf}%
}%
& 
{\includegraphics[
height=0.3269in,
width=0.3269in
]%
{ut10.pdf}%
}%
\end{array}
$}}
\quad\underset{K_{95}=\text{$0_1+0_1+0_1$}}{\fbox{$%
\begin{array}
[c]{ccc}
{\includegraphics[
height=0.3269in,
width=0.3269in
]%
{ut09.pdf}%
}%
& {\includegraphics[
height=0.3269in,
width=0.3269in
]%
{ut09.pdf}%
}%
 \\
{\includegraphics[
height=0.3269in,
width=0.3269in
]%
{ut07.pdf}%
}%
& 
{\includegraphics[
height=0.3269in,
width=0.3269in
]%
{ut09.pdf}%
}%
\end{array}
$}}
\quad\underset{K_{96}=\text{$0_1+0_1+0_1$}}{\fbox{$%
\begin{array}
[c]{ccc}
{\includegraphics[
height=0.3269in,
width=0.3269in
]%
{ut09.pdf}%
}%
& {\includegraphics[
height=0.3269in,
width=0.3269in
]%
{ut09.pdf}%
}%
 \\
{\includegraphics[
height=0.3269in,
width=0.3269in
]%
{ut08.pdf}%
}%
& 
{\includegraphics[
height=0.3269in,
width=0.3269in
]%
{ut09.pdf}%
}%
\end{array}
$}}
\\ \\ \\
$\bigskip$
\hspace{-0.1in}\underset{K_{97}=\text{$0_1+0_1+0_1+0_1$}}{\fbox{$%
\begin{array}
[c]{ccc}
{\includegraphics[
height=0.3269in,
width=0.3269in
]%
{ut09.pdf}%
}%
& {\includegraphics[
height=0.3269in,
width=0.3269in
]%
{ut09.pdf}%
}%
 \\
{\includegraphics[
height=0.3269in,
width=0.3269in
]%
{ut09.pdf}%
}%
& 
{\includegraphics[
height=0.3269in,
width=0.3269in
]%
{ut09.pdf}%
}%
\end{array}
$}}
\ \; \underset{K_{98}=\text{$8_3^4$ (L8n8)}}{\fbox{$%
\begin{array}
[c]{ccc}
{\includegraphics[
height=0.3269in,
width=0.3269in
]%
{ut10.pdf}%
}%
& {\includegraphics[
height=0.3269in,
width=0.3269in
]%
{ut10.pdf}%
}%
 \\
{\includegraphics[
height=0.3269in,
width=0.3269in
]%
{ut10.pdf}%
}%
& 
{\includegraphics[
height=0.3269in,
width=0.3269in
]%
{ut10.pdf}%
}%
\end{array}
$}}
\quad\underset{K_{99}=\text{$2_1^2\#2_1^2+0_1$}}{\fbox{$%
\begin{array}
[c]{ccc}
{\includegraphics[
height=0.3269in,
width=0.3269in
]%
{ut10.pdf}%
}%
& {\includegraphics[
height=0.3269in,
width=0.3269in
]%
{ut09.pdf}%
}%
 \\
{\includegraphics[
height=0.3269in,
width=0.3269in
]%
{ut10.pdf}%
}%
& 
{\includegraphics[
height=0.3269in,
width=0.3269in
]%
{ut09.pdf}%
}%
\end{array}
$}}
\quad\underset{K_{100}=\text{$2_1^2+0_1+0_1$}}{\fbox{$%
\begin{array}
[c]{ccc}
{\includegraphics[
height=0.3269in,
width=0.3269in
]%
{ut09.pdf}%
}%
& {\includegraphics[
height=0.3269in,
width=0.3269in
]%
{ut10.pdf}%
}%
 \\
{\includegraphics[
height=0.3269in,
width=0.3269in
]%
{ut09.pdf}%
}%
& 
{\includegraphics[
height=0.3269in,
width=0.3269in
]%
{ut09.pdf}%
}%
\end{array}
$}}
\quad\underset{K_{101}=\text{$2_1^2\#2_1^2+0_1$}}{\fbox{$%
\begin{array}
[c]{ccc}
{\includegraphics[
height=0.3269in,
width=0.3269in
]%
{ut10.pdf}%
}%
& {\includegraphics[
height=0.3269in,
width=0.3269in
]%
{ut10.pdf}%
}%
 \\
{\includegraphics[
height=0.3269in,
width=0.3269in
]%
{ut09.pdf}%
}%
& 
{\includegraphics[
height=0.3269in,
width=0.3269in
]%
{ut09.pdf}%
}%
\end{array}
$}}
\\ \\ \\
$\bigskip$
\underset{K_{102}=\text{$2_1^2\#2_1^2\#2_1^2$}}{\fbox{$%
\begin{array}
[c]{ccc}
{\includegraphics[
height=0.3269in,
width=0.3269in
]%
{ut10.pdf}%
}%
& {\includegraphics[
height=0.3269in,
width=0.3269in
]%
{ut10.pdf}%
}%
 \\
{\includegraphics[
height=0.3269in,
width=0.3269in
]%
{ut09.pdf}%
}%
& 
{\includegraphics[
height=0.3269in,
width=0.3269in
]%
{ut10.pdf}%
}%
\end{array}
$}}
\quad\underset{K_{103}=\text{ $2_1^2+2_1^2$}}{\fbox{$%
\begin{array}
[c]{ccc}
{\includegraphics[
height=0.3269in,
width=0.3269in
]%
{ut10.pdf}%
}%
& {\includegraphics[
height=0.3269in,
width=0.3269in
]%
{ut09.pdf}%
}%
 \\
{\includegraphics[
height=0.3269in,
width=0.3269in
]%
{ut09.pdf}%
}%
& 
{\includegraphics[
height=0.3269in,
width=0.3269in
]%
{ut10.pdf}%
}%
\end{array}$}}
\quad\underset{K_{104}=\text{$2_1^2\#2_1^2+0_1$}}{\fbox{$%
\begin{array}
[c]{ccc}
{\includegraphics[
height=0.3269in,
width=0.3269in
]%
{ut10.pdf}%
}%
& {\includegraphics[
height=0.3269in,
width=0.3269in
]%
{ut09.pdf}%
}%
 \\
{\includegraphics[
height=0.3269in,
width=0.3269in
]%
{ut10.pdf}%
}%
& 
{\includegraphics[
height=0.3269in,
width=0.3269in
]%
{ut09.pdf}%
}%
\end{array}
$}}
$

\end{document}